\documentclass[a4paper,11pt]{amsart}
\usepackage{amsfonts,stmaryrd,amssymb,amsmath,a4wide,verbatim,hyperref,url,color,epsfig,mathrsfs,graphicx,pgf,tikz,amsthm,amsfonts,mathtext,cite,enumerate,float}
\usepackage[british]{babel}

\usepackage[active]{srcltx}
\usepackage[T1,T2A]{fontenc}
\usepackage[utf8]{inputenc}

\usetikzlibrary{arrows}
\usepackage[matrix,arrow,curve]{xy}
\usetikzlibrary{arrows}

\makeatletter
\@addtoreset{equation}{section}\makeatother

\theoremstyle{plain}
\newtheorem{theorem}{Theorem}[section]

\theoremstyle{definition}
\newtheorem{definition}{Definition}[section]
\theoremstyle{definition}
\newtheorem{commentary}{Remark}[section]

\begin{document}

\title[Some new theorems in plane geometry]{Some new theorems in plane geometry}

\subjclass[2010]{51M04, 51N20}

\keywords{Plane geometry, Analytic geometry}

\author[Alexander Skutin]{Alexander Skutin}

\maketitle

\section{Introduction}

In this article we will represent some ideas and a lot of new theorems in plane geometry.

\bigskip

\section{Deformation of equilateral triangle}

\subsection{Deformation principle for equilateral triangle} If some triangle points lie on a circle (line) or are equivalent in the case of an equilateral triangle, then in the general case of an arbitrary triangle they are connected by some natural relations. Thus, triangle geometry can be seen as a deformation of the equilateral triangle geometry.

This principle partially describes why there exists so many relations between the Kimberling centers $X_i$ (see definition of the Kimberling center here \cite{Kim}).

\subsection{First example of application} Consider an equilateral triangle $ABC$ with center $O$ and let $P$ be an arbitrary point on the circle $(ABC)$. Then we can imagine that $O$ is the first Fermat point of the triangle $ ABC $, and the point $P$ is the second Fermat point of the triangle $ABC$. Consider a circle $\omega$ with center at $P$ and passing through the point $O$. Let $X$, $Y$ be the intersection points of the circles $(ABC)$ and $\omega$. Then it is easy to see that the line $XY$ passes through the middle of the segment $OP$. Thus, we can guess that this fact takes place in general for an arbitrary triangle $ABC$, and it's first and the second Fermat points $F_1$, $F_2$. Thus, we obtain the following theorem.

\begin{theorem}

Consider triangle $ABC$ with the first Fermat point $F_1$ and the second Fermat point $F_2$. Consider a circle $\omega$ with center at $F_2$ and radius $F_2F_1$. Then the radical line of circles $\omega$, $(ABC)$ goes through the midpoint of the segment $F_1F_2$.

\end{theorem}

\definecolor{ffffff}{rgb}{1.0,1.0,1.0}
\definecolor{ffwwqq}{rgb}{1.0,0.4,0.0}
\definecolor{ffqqtt}{rgb}{1.0,0.0,0.2}
\begin{tikzpicture}[line cap=round,line join=round,>=triangle 45,x=1.0cm,y=1.0cm]
%\clip(-4.3,-2.2800000000000002) rectangle (11.16,6.3);
\draw(-1.0326794919243112,1.1857219549304678) circle (1.951751350283466cm);
\draw(-1.9456925538575303,2.9107554276702468) circle (1.9517513502834662cm);
\draw [dash pattern=on 3pt off 3pt] (-2.98310883266206,1.2575461856791745)-- (0.004736786880218679,2.83893119692154);
\draw (-1.9456925538575303,2.9107554276702468)-- (-1.0326794919243112,1.1857219549304678);
\draw (-2.74,0.24)-- (-0.9980384757729328,3.1371658647914034);
\draw (-0.9980384757729328,3.1371658647914034)-- (0.64,0.18);
\draw (-2.74,0.24)-- (0.64,0.18);
\draw (4.06,0.54)-- (5.78,3.38);
\draw (5.78,3.38)-- (7.58,0.28);
\draw (4.06,0.54)-- (7.58,0.28);
\draw(5.891049406357717,1.3718996553044809) circle (2.01116855708781cm);
\draw(4.6638432799308625,3.0499539049055175) circle (1.9325894814553557cm);
\draw [dash pattern=on 3pt off 3pt] (3.8818627124472997,1.2826371546001398)-- (6.585034591306983,3.2595393264519172);
\draw (4.6638432799308625,3.0499539049055175)-- (5.720051721346142,1.4315207638351217);
\begin{scriptsize}
\draw [fill=ffqqtt] (-2.74,0.24) circle (1.5pt);
\draw [fill=ffqqtt] (0.64,0.18) circle (1.5pt);
\draw [fill=ffqqtt] (-0.9980384757729328,3.1371658647914034) circle (1.5pt);
\draw [fill=ffwwqq] (-1.9456925538575303,2.9107554276702468) circle (1.5pt);
\draw[color=ffwwqq] (-1.8,3.16) node {$P$};
\draw [fill=ffwwqq] (-1.0326794919243112,1.1857219549304678) circle (1.5pt);
\draw[color=ffwwqq] (-0.9199999999999999,1.1199999999999999) node {$O$};
\draw [fill=ffffff] (-2.98310883266206,1.2575461856791745) circle (1.5pt);
\draw [fill=ffffff] (0.004736786880218679,2.83893119692154) circle (1.5pt);
\draw [fill=ffqqtt] (4.06,0.54) circle (1.5pt);
\draw [fill=ffqqtt] (5.78,3.38) circle (1.5pt);
\draw [fill=ffqqtt] (7.58,0.28) circle (1.5pt);
\draw [fill=ffwwqq] (5.720051721346142,1.4315207638351217) circle (1.5pt);
\draw[color=ffwwqq] (5.96,1.44) node {$F_1$};
\draw [fill=ffwwqq] (4.6638432799308625,3.0499539049055175) circle (1.5pt);
\draw[color=ffwwqq] (4.84,3.3400000000000003) node {$F_2$};
\draw [fill=ffffff] (3.8818627124472997,1.2826371546001398) circle (1.5pt);
\draw [fill=ffffff] (6.585034591306983,3.2595393264519172) circle (1.5pt);
\draw [fill=ffffff] (5.1919475006385,2.240737334370321) circle (1.5pt);
\draw [fill=ffffff] (-1.4891860228909208,2.0482386913003574) circle (1.5pt);
\end{scriptsize}
\end{tikzpicture}

\subsection{Fermat points}

\begin{theorem}

Consider triangle $ABC$ with the first Fermat point $F_1$ and the second Fermat point $F_2$. Consider a point $A'$ on $BC$, such that $F_2A'$ is parallel to $AF_1$. Similarly define $B'$, $C'$. Then $A'$, $B'$, $C'$ lie on the same line which goes through the middle of $F_1F_2$.

\end{theorem}

\begin{theorem}

Given triangle $ABC$ with the first Fermat point $F_1$ and the second Fermat point $F_2$. Let $F_A$ be the second Fermat point of $BCF_1$, similarly define $F_B$, $F_C$.

\begin{enumerate}

\item Reflect $F_A$ wrt line $AF_1$ and get the point $F_A'$, similarly define the points $F_B'$, $F_C'$. Then $F_2$ lies on the circles $(F_AF_BF_C)$, $(F_A'F_B'F_C')$.

\item Consider isogonal conjugations $F_A^A$, $F_A^B$, $F_A^C$ of $F_A$ wrt $BCF_1$, $ACF_1$, $ABF_1$ respectively. Then $\angle F_A^BF_A^AF_A^C = \pi/3$.

\item Triangles $ABC$ and $F_A^AF_B^BF_C^C$ are perspective.

\end{enumerate}

\end{theorem}

\begin{theorem}

Let $F_1$, $F_2$ be the first and the second Fermat points for $ABC$. Consider an alphabet with the letters $"a, b, c"$. Consider the set $\mathcal{W}$ of all words in this alphabet, and let the function $|\cdot|$ denote the length of a word in $\mathcal{W}$. Inductively define $F_A^{(\omega)}$, $F_B^{(\omega)}$, $F_C^{(\omega)}$ :

\begin{enumerate}[(i)]

\item Let $F_A^{(\emptyset)}= F_1$, similarly for $F_B^{(\emptyset)}$ and $F_C^{(\emptyset)}$.

\item Let $F_A^{(\omega a)}$ be the second Fermat point of $BCF_A^{(\omega)}$, same for $F_B^{(\omega a)}$, $F_C^{(\omega a)}$.

\item By definition, let $F_A^{(\omega b)}$ be a reflection of the point $F_A^{(\omega)} $ wrt side $BC$. Similarly define the points $F_B^{(\omega b)}$, $F_C^{(\omega b)}$.

\item By definition, let $F_A^{(\omega c)}$ be a reflection of the point $F_A^{(\omega)}$ wrt midpoint of $BC$. Similarly define the points $F_B^{(\omega c)}$, $F_C^{(\omega c)}$.

\end{enumerate}

Then we have that :

\begin{enumerate}

\item If $|\omega|$ is odd, then the points $F_2$, $F_A^{(\omega)}$, $F_B^{(\omega)}$, $F_C^{(\omega)}$ lie on the same circle.

\item If $|\omega|$ is even, then the points $F_1$, $F_A^{(\omega)}$, $F_B^{(\omega)}$, $F_C^{(\omega)}$ lie on the same circle.

\end{enumerate}

\end{theorem}

\begin{theorem}

Given triangle $ABC$ with the first Fermat point $F$. Let $F_A$ be the second Fermat point for $FBC$. Similarly define $F_B$, $F_C$. Let $F_B^{A}$ be the first Fermat point of $FF_BA$ and $F_C^{A}$ be the first Fermat point of $FF_CA$. Let two tangents from $F_B^{A}$, $F_C^{A}$ to circle $(F_B^{A}F_C^{A}A)$ meet at $X_A$. Similarly define $X_B$, $X_C$. Then $F$, $X_A$, $X_B$, $X_C$ are collinear.

\end{theorem}

\subsection{Miquel points}

\begin{theorem}

Consider any triangle $ABC$ and any point $X$. Let $P= AX\cap BC$, $Q = BX\cap AC$, $R = CX\cap AB$. Let $M_Q$ be the Miquel point of the lines $XA$, $XC$, $BA$, $BC$, similarly define the points $M_P$, $M_R$. Denote $l_Q$ as the Simson line of the point $M_Q$ wrt triangle $XPC$ (it's same for the triangles $XRA$, $APB$, $RBC$). Similarly define the lines $l_P$, $l_R$. Let $P'Q'R'$ be the midpoint triangle for $PQR$. Let $P'Q'\cap l_R = R^*$, like the same define the points $P^*, Q^*$. Then $P'P^*$, $Q'Q^*$, $R'R^*$ are concurrent.

\end{theorem}

\subsection{Steiner lines}

\begin{theorem}

Consider any triangle $ABC$ and any point $P$. Let $\mathcal{L}
_{A}$ be the Steiner line for the lines $AB$, $AC$, $PB$, $PC$, similarly define the lines $\mathcal{L}_{B}$, $\mathcal{L}_{C}$. Let the lines $\mathcal{L}_{A}$, $\mathcal{L}_{B}$, $\mathcal{L}_{C}$ form a triangle $\triangle$. Then the circumcircle of $\triangle$, the pedal circle of $P$ wrt $ABC$ and the nine~- point circle of $ABC$ passes through the same point.

\end{theorem}

\subsection{Isogonal conjugations}

\begin{theorem}

Let the incircle of $ABC$ touches sides $BC$, $CA$, $AB$ at $A_1$, $B_1$, $C_1$ respectively. Let the $A$, $B$, $C$~-- excircles of $ABC$ touches sides $BC$, $CA$, $AB$ at $A_2$, $B_2$, $C_2$ respectively. Let $A^{1}$ be isogonal conjugate to the point $A$ wrt $A_1B_1C_1$. Similarly define the points $B^{1}$ and $C^{1}$. Let $A^{2}$ be isogonal conjugate to the point $A$ wrt $A_2B_2C_2$. Similarly define the points $B^{2}$ and $C^{2}$. Then we have that :

\begin{enumerate}

\item Lines $A_2A^{1}$, $B_2B^{1}$, $C_2C^{1}$ are concurrent at $X_2^{1}$.

\item Lines $A_1A^{2}$, $B_1B^{2}$, $C_1C^{2}$ are concurrent at $X_1^{2}$.

\item Points $A^{1}$, $B^{1}$, $C^{1}$, $A^{2}$, $B^{2}$, $C^{2}$, $X_2^{1}$, $X_1^{2}$ lie on the same conic.

\end{enumerate}

\end{theorem}

\subsection{Deformation principle for circles} If in the particular case of an equilateral triangle some two circles are equivalent, then in the general case the radical line of these two circles has some natural relations with respect to the base triangle.

\begin{theorem}

Let given triangle $ABC$ and its centroid $G$. Let the circumcenters of the triangles $ABG$, $CBG$, $AGC$ form a triangle with the circumcircle $\omega$. Then the circumcenter of the pedal triangle of $G$ wrt $ABC$ lies on the radical line of $(ABC)$ and $\omega$.

\end{theorem}

\begin{theorem}

Consider triangle $ABC$ with the first Fermat point $F$ and let $O_A$, $O_B$, $O_C$ be the circumcenters of $FBC$, $FAC$, $FAB$ respectively. Let $X_A$ be a reflection of $O_A$ wrt $BC$, similarly define the points $X_B$, $X_C$. Let $l$ be the radical line of $(O_AO_BO_C)$ and $(ABC)$. Then $F$ is the Miquel point of the lines $X_AX_B$, $X_BX_C$, $X_AX_C$, $l$.

\end{theorem}
\begin{commentary}
In all theorems from this section consider the case of an equilateral triangle.

\end{commentary}
\bigskip
\section{Nice fact about triangle centers}

For any $i, j\in\mathbb{N}$ if some fact is true for the Kimberling center $X_i$, then we can try to transport the same construction from the triangle center $X_i$ to the center $X_j$ and to look on the nice properties of the resulting configuration.

\begin{theorem}

Given a triangle $ABC$. Let $A$~-- excircle tangent to $BC$ at $A_1$. Similarly define $B_1$, $C_1$. Let $A_2B_2C_2$ be the midpoint triangle of $ABC$. Consider intersection $N$ of the perpendiculars to $AB$, $BC$, $CA$ from $C_1$, $A_1$, $B_1$ respectively. Let the circle $(A_1B_1C_1)$ meet sides $AB$, $BC$ second time at $C_3$, $A_3$ respectively. Let $X= A_3C_3\cap B_2C_2$. Then $AX\perp CN$.

\end{theorem}

\begin{commentary}

Consider the case when $N$ is the incenter of $ABC$.

\end{commentary}

\begin{theorem}

Let the incircle of $ABC$ touches sides $AB$, $BC$, $CA$ at $C'$, $A'$, $B'$ respectively. Reflect $B$, $C$ wrt line $AA'$ and get the points $B_1$, $C_1$ respectively. Let $A^*= CB_1\cap BC_1$, similarly define the points $B^*$, $C^*$. Then the points $A^*$, $B^*$, $C^*$, $I$, $G$ lie on the same circle, where $G$ is the Gergonne point of $ABC$.

\end{theorem}

\begin{commentary}

Consider the case when $AA'$, $BB'$, $CC'$ goes through the first Fermat point.

\end{commentary}

\begin{theorem}

Let the incircle of $ABC$ touches sides $AB$, $BC$, $CA$ at $C'$, $A'$, $B'$ respectively. Reflect point $A'$ wrt lines $BB'$, $CC'$ and get the points $A_B$, $A_C$ respectively. Let $A^*= A_BC'\cap A_CB'$, similarly define the points $B^*$, $C^*$. Then the incenter of $ABC$ coincide with the circumcenter of $A^*B^*C^*$.

\end{theorem}

\section{Construction of midpoint analog}

\begin{definition}

For any pairs of points $A$, $B$ and $C$, $D$ denote $\mathcal{M}(AB, CD)$ as the Miquel point of the complete quadrilateral formed by the four lines $AC$, $AD$, $BC$, $BD$.

\end{definition}

\begin{definition}

For any point $X$ and a segment $YZ$ denote $\mathcal{M}(X, YZ)$~-- as a point $P$, such that the circles $(PXY)$ and $(PXZ)$ are tangent to segments $XZ$, $XY$ at $X$.

\end{definition}

Consider any two segments $AB$ and $CD$, then the point $\mathcal{M}(AB, CD)$ can be seen as midpoint between the two segments $AB$, $CD$. Also we can consider the segment $AB$ and the point $C$ and to look on the point $\mathcal{M}(C, AB)$ as on the midpoint between the point $C$ and the segment $AB$.

\begin{commentary}

In the case when $A=B$ and $C=D$ we will get that the points $\mathcal{M}(C, AB)$ and $\mathcal{M}(AB, CD)$ are midpoints of $AC$.

\end{commentary}

\begin{theorem}

Let given segments $P_AQ_A$, $P_BQ_B$, $P_CQ_C$, $P_DQ_D$. Let $\omega_D$ be the circumcircle of $\mathcal{M}(P_AQ_A, P_BQ_B)\mathcal{M}(P_AQ_A, P_CQ_C)\mathcal{M}(P_CQ_C, P_BQ_B)$. Like the same define the circles $\omega_A$, $\omega_B$, $\omega_C$. Then the circles $\omega_A$, $\omega_B$, $\omega_C$, $\omega_D$ intersect at the same point.

\end{theorem}

\begin{theorem}[\textbf{Three nine-point circles intersect}]

Let given circle $\omega$ and points $P_A$, $Q_A$, $P_B$, $Q_B$, $P_C$, $Q_C$ on it. Let $A'= \mathcal{M}(P_BQ_B, P_CQ_C)$, similarly define points $B'$, $C'$. Consider points $A_1= P_BP_C\cap Q_BQ_C$, $B_1 = P_AP_C\cap Q_AQ_C$, $C_1 = P_AP_B\cap Q_AQ_B$. Let the lines $A_1A'$, $B_1B'$, $C_1C'$ form a triangle with the nine~- point circle $\omega_1$. Then the circles $\omega_1$, $(A'B'C')$, $(A_1B_1C_1)$ pass through the same point.

\end{theorem}

\begin{theorem}

Consider triangle $ABC$. Let points $P$, $Q$ considered such that $PP'\parallel QQ'$, where $P$, $P'$ are isogonal conjugated and $Q$, $Q'$ are isogonal conjugated  wrt triangle $ABC$. Let $P_A$ be the second intersection point of line $AP'$ with circle $(BCP')$. Similarly define the points $P_B$, $P_C$. Let $Q_A$ be the second intersection point of line $AQ'$ with circle $(BCQ')$. Similarly define the points $Q_B$, $Q_C$. Let $Q_AP_B\cap Q_BP_A = R_C$, similarly $Q_AP_C\cap Q_CP_A = R_B$ and $Q_BP_C\cap Q_CP_B = R_A$. Then we have that :

\begin{enumerate}

\item Circles $(\mathcal{M}(Q_AP_A , Q_BP_B)\mathcal{M}(Q_CP_C , Q_BP_B)\mathcal{M}(Q_AP_A , Q_CP_C))$ $$(\mathcal{M}(Q_AP_A , Q_BP_B)\mathcal{M}(Q_AP_A , QP)\mathcal{M}(QP , Q_BP_B))$$ and $(R_AR_BR_C)$ goes through the same point.

\item Lines $AR_A$, $BR_B$, $CR_C$ intersect at the same point which lies on the circle $$(\mathcal{M}(Q_AP_A , QP)\mathcal{M}(Q_BP_B , QP)\mathcal{M}(QP , Q_CP_C))$$

\end{enumerate}
\end{theorem}
\begin{commentary}

Consider the case when $P = Q$~-- incenter of $ABC$.

\end{commentary}

\begin{theorem}

Let $P$ and $Q$ be two isogonal conjugated points wrt $ABC$. Let $O_A^P$, $O_B^P$, $O_C^P$ be the circumcenters of the triangles $BCP$, $CAP$, $ABP$ respectively. Similarly let $O_A^Q$, $O_B^Q$, $O_C^Q$ be the circumcenters of the triangles $BCQ$, $CAQ$, $ABQ$ respectively. Let $M_A = \mathcal{M}(O_B^PO_B^Q , O_C^PO_C^Q)$, similarly define the points $M_B$, $M_C$. Let $O$ be the circumcenter of $ABC$. Then we have that :

\begin{enumerate}

\item Circles $(OM_AA)$, $(OM_BB)$, $(OM_CC)$ are coaxial.

\item Let $A^{\triangle}=\mathcal{M}(PQ , O_A^PO_A^Q)$, $B^{\triangle}=\mathcal{M}(PQ , O_B^PO_B^Q)$ and $C^{\triangle}=\mathcal{M}(PQ , O_C^PO_C^Q)$. Then $A^{\triangle}$ lies on $(OM_AA)$.

\item Points $A^{\triangle}$, $B^{\triangle}$, $C^{\triangle}$, $M_A$, $M_B$, $M_C$ lie on the same circle.

\end{enumerate}

\end{theorem}

\section{Deformation of segment into a conic}

\subsection{Deformation principle for conics} In some statements we can replace some segment by a conic.

\begin{theorem}

Consider a circle $\omega$ and the two conics $\mathcal{K}_1$, $\mathcal{K}_2$ which are tangent to $\omega$ at four points $P_1$, $Q_1$, $P_2$, $Q_2$. Where $P_1$ and $Q_1$ lie on $\mathcal{K}_1$, $P_2$ and $Q_2$ lie on $\mathcal{K}_2$. Consider two external tangents to $\mathcal{K}_1$, $\mathcal{K}_2$ which intersect at the point $C$. Let the circles $(CP_1Q_1)$ and $(CP_2Q_2)$ intersect second time at $E$. Then $\angle CEO =\pi/2$, where $O$ is the center of $\omega$.

\end{theorem}

\begin{theorem}

Consider a circle $\omega$ and the two conics $\mathcal{K}_1$, $\mathcal{K}_2$ which are tangent to $\omega$ at four points. Let $\mathcal{K}_1$ has foci $F_1$, $F_2$ and $\mathcal{K}_2$ has foci $F_3$, $F_4$. Let $l_1$, $l_2$ be two external tangents to $\mathcal{K}_1$, $\mathcal{K}_2$. Let $Y= F_1F_4\cap F_2F_3$ and the circles $(F_1YF_3)$, $(F_2YF_4)$ meet second time at $X$. Consider the point $Z= l_1\cap l_2$ (see picture below). Then $X$, $Y$, $Z$ are collinear.

\end{theorem}
\medskip
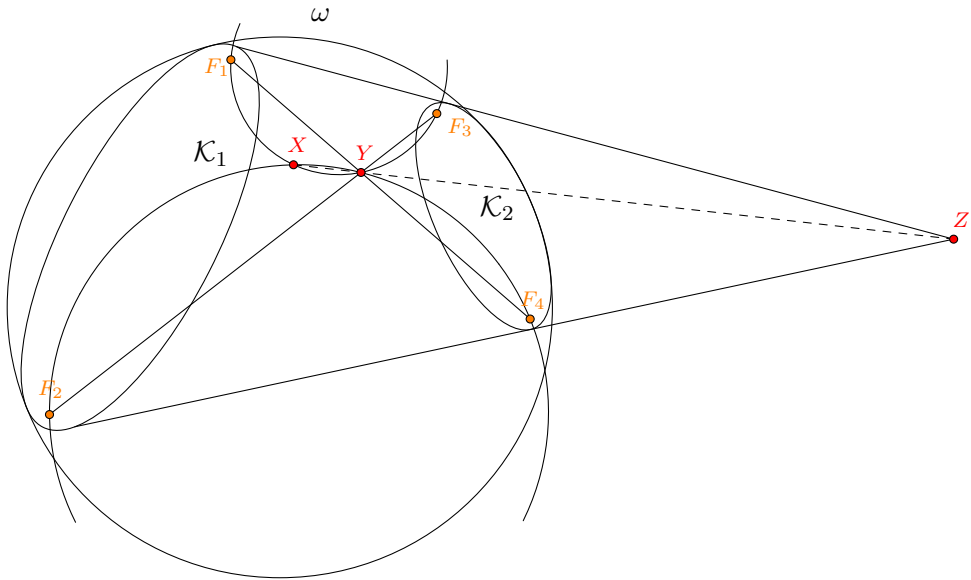
\definecolor{ffqqqq}{rgb}{1.0,0.0,0.0}
\definecolor{ffxfqq}{rgb}{1.0,0.4980392156862745,0.0}
\begin{tikzpicture}[line cap=round,line join=round,>=triangle 45,x=1.0cm,y=1.0cm]
%\clip(-5.261634057697279,-7.19310444030945) rectangle (11.215710178777087,2.9004068426689793);
\draw(0.21764574025017722,-3.4143179245217703) circle (3.5851591753888847cm);
\draw [rotate around={63.110919658012214:(-1.6184878332753352,-2.483234159309076)}] (-1.6184878332753352,-2.483234159309076) ellipse (2.826149470384081cm and 1.015214678769174cm);
\draw [rotate around={-65.73472634378953:(2.896376658748347,-2.206778542255469)}] (2.896376658748347,-2.206778542255469) ellipse (1.6280316690209338cm and 0.6468347049734543cm);
\draw (-1.051262947548442,-1.0823766399042454) node[anchor=north west] {$\mathcal{K}_1$};
\draw (2.71531228882795,-1.7992776667674257) node[anchor=north west] {$\mathcal{K}_2$};
\draw (0.496332919965725,0.6814274738067541) node[anchor=north west] {$\omega$};
\draw (-2.8113406888643215,-4.835586609491031)-- (2.2823918935476546,-0.8447525924738806);
\draw (3.5103614239490377,-3.5688044920370543)-- (-0.4256349776863491,-0.1308817091271236);
\draw (9.079635740061052,-2.509921792687575)-- (-0.35390814946863464,0.04390349832076328);
\draw (-2.4796096589712375,-4.998127766278186)-- (9.079635740061052,-2.509921792687575);
\draw [dash pattern=on 3pt off 3pt] (0.3977068654526339,-1.5235055057251057)-- (9.079635740061052,-2.509921792687575);
\draw [shift={(0.4701783844458631,-4.804372710432985)}] plot[domain=2.6429194168184713:3.6029264988446057,variable=\t]({1.0*3.281667523682577*cos(\t r)+-0.0*3.281667523682577*sin(\t r)},{0.0*3.281667523682577*cos(\t r)+1.0*3.281667523682577*sin(\t r)});
\draw [shift={(0.47017838444584836,-4.804372710432978)}] plot[domain=1.5928818703158307:2.642919416818471,variable=\t]({1.0*3.2816675236825605*cos(\t r)+-0.0*3.2816675236825605*sin(\t r)},{0.0*3.2816675236825605*cos(\t r)+1.0*3.2816675236825605*sin(\t r)});
\draw [shift={(0.4701783844458578,-4.804372710432995)}] plot[domain=0.6415567997601213:1.5928818703158334,variable=\t]({1.0*3.281667523682578*cos(\t r)+-0.0*3.281667523682578*sin(\t r)},{0.0*3.281667523682578*cos(\t r)+1.0*3.281667523682578*sin(\t r)});
\draw [shift={(0.4701783844458389,-4.804372710432986)}] plot[domain=-0.45416920282476436:0.6415567997601157,variable=\t]({1.0*3.281667523682587*cos(\t r)+-0.0*3.281667523682587*sin(\t r)},{0.0*3.281667523682587*cos(\t r)+1.0*3.281667523682587*sin(\t r)});
\draw [shift={(0.9963786156781685,-0.22986258405152316)}] plot[domain=2.723428314916651:4.9167687559035045,variable=\t]({1.0*1.425454269106613*cos(\t r)+-0.0*1.425454269106613*sin(\t r)},{0.0*1.425454269106613*cos(\t r)+1.0*1.425454269106613*sin(\t r)});
\draw [shift={(0.996378615678169,-0.2298625840515243)}] plot[domain=-1.3664165512760817:0.06765800131192408,variable=\t]({1.0*1.4254542691066117*cos(\t r)+-0.0*1.4254542691066117*sin(\t r)},{0.0*1.4254542691066117*cos(\t r)+1.0*1.4254542691066117*sin(\t r)});
\begin{scriptsize}
\draw [fill=ffxfqq] (-0.4256349776863491,-0.1308817091271236) circle (1.5pt);
\draw[color=ffxfqq] (-0.6074670737759971,-0.2289230364956972) node {$F_1$};
\draw [fill=ffxfqq] (-2.8113406888643215,-4.835586609491031) circle (1.5pt);
\draw[color=ffxfqq] (-2.79230829850188,-4.496191053538438) node {$F_2$};
\draw [fill=ffxfqq] (3.5103614239490377,-3.5688044920370543) circle (1.5pt);
\draw[color=ffxfqq] (3.5460071294789364,-3.324114771524032) node {$F_4$};
\draw [fill=ffxfqq] (2.2823918935476546,-0.8447525924738806) circle (1.5pt);
\draw[color=ffxfqq] (2.590139093661363,-1.025479733010342) node {$F_3$};
\draw [fill=ffqqqq] (1.2856886476409826,-1.6256489497679758) circle (1.5pt);
\draw[color=ffqqqq] (1.3384071419954924,-1.3668611743737613) node {$Y$};
\draw [fill=ffqqqq] (0.3977068654526339,-1.5235055057251057) circle (1.5pt);
\draw[color=ffqqqq] (0.4508153944506024,-1.2644467419647356) node {$X$};
\draw [fill=ffqqqq] (9.079635740061052,-2.509921792687575) circle (1.5pt);
\draw[color=ffqqqq] (9.178800911975353,-2.2544529219186513) node {$Z$};
\end{scriptsize}
\end{tikzpicture}
\medskip
\begin{definition}
For any two conics $\mathcal{C}_1$ and $\mathcal{C}_2$ consider all their four tangents $l_1$, $l_2$, $l_3$, $l_4$. By definition let $\mathcal{M}(\mathcal{C}_1 , \mathcal{C}_2)$ be the Miquel point for the lines $l_1$, $l_2$, $l_3$, $l_4$.

\end{definition}

\begin{theorem}

Let given conics $\mathcal{C}_A$, $\mathcal{C}_B$, $\mathcal{C}_C$, $\mathcal{C}_D$. Let $\omega_D$ be the circumcircle of triangle $$\mathcal{M}(\mathcal{C}_A, \mathcal{C}_B)\mathcal{M}(\mathcal{C}_A, \mathcal{C}_C)\mathcal{M}(\mathcal{C}_C, \mathcal{C}_B)$$ Like the same define $\omega_A$, $\omega_B$, $\omega_C$. Then $\omega_A$, $\omega_B$, $\omega_C$, $\omega_D$ intersect at the same point.

\end{theorem}

\begin{theorem}

Consider any triangle $ABC$. Let $P$, $P'$ and $Q$, $Q'$ be two pairs of isogonal conjugated points wrt $ABC$. Consider a conic $\mathcal{C}_P$ with foci at $P$, $P'$ and also consider a conic $\mathcal{C}_Q$ with foci at $Q$, $Q'$. Then the points $A$, $B$, $C$, $\mathcal{M}(\mathcal{C}_P, \mathcal{C}_Q)$ lie on the same circle.

\end{theorem}

\newpage

\section{On incircles in orthocenter construction}

Main idea of this section is to construct something nice which includes incircles in orthocenter construction.

\begin{theorem}

Consider triangle $ABC$ with orthocenter $H$ and altitudes $AA'$, $BB'$, $CC'$. Let $C^*=A'B'\cap CC'$. Let $I_1$, $I_2$, $I_3$, $I_4$, $I_5$, $I_6$ be the incenters of the triangles $AB'C'$, $BA'C'$, $HBA'$, $HAB'$, $HA'C^*$, $HC^*C'$ respectively. Then we have that :

\begin{enumerate}

\item Lines $I_1I_4$, $I_2I_3$, $CC'$ are concurrent.

\item Lines $I_6I_4$, $I_5I_3$, $CC'$ are concurrent.

\end{enumerate}

\end{theorem}

\begin{theorem}

Let given triangle $ABC$ with orthocenter $H$ and altitudes $AH_A$, $BH_B$, $CH_C$. Let the second tangent line through $H_A$ to $A$~-- excircle of $ABC$ meet the second tangent line through $H_B$ to $B$~-- excircle of $ABC$ at $C_1$. Let the second tangent line through $H_A$ to $C$~-- excircle of $ABC$ meet the second tangent line through $H_B$ to $C$~-- excircle of $ABC$ at $C_2$. Then $H_C$, $C_1$, $C_2$ are collinear.

\end{theorem}

\begin{theorem}

Let given triangle $ABC$ with altitudes $AH_A$, $BH_B$, $CH_C$. Let $l_A$ be the line which passes through the incenters of $AH_AB$ and $AH_AC$. Similarly define $l_B$, $l_C$. Then the triangle formed by $l_A$, $l_B$, $l_C$ is perspective to $ABC$.

\end{theorem}

\begin{theorem}

Let $I$ be the incenter of $ABC$. Let $H_a$ be the orthocenter of $BCI$, similarly define $H_b$, $H_c$. Let $I_a$ be the incenter of $BH_aC$, similarly define $I_b$ and $I_c$. Then

\begin{enumerate}

\item Point $I$ is the orthocenter of $I_aI_bI_c$.

\item Let $H_a^{(2)}$ be the orthocenter of $BI_aC$, similarly define $H_b^{(2)}$, $H_c^{(2)}$. Then the triangles $H_a^{(2)}H_b^{(2)}H_c^{(2)}$ and $ABC$ are perspective.

\item Let $I_a^{(2)}$ be the incenter of $BH_a^{(2)}C$, similarly define $I_b^{(2)}$, $I_c^{(2)}$. Then the triangles $ABC$ and $I_a^{(2)}I_b^{(2)}I_c^{(2)}$ are perspective.

\end{enumerate}

\end{theorem}

\begin{theorem}[\textbf{Four nine-point circles intersect}]

Let given triangle $ABC$ with orthocenter $H$. Let $\omega$ be the incircle of $ABC$ and $\omega_A$ be the incircle of $BCH$. Like the same define the circles $\omega_B$, $\omega_C$. Let $l_A$ be the second external tangent line of the circles $\omega_A$, $\omega$. Similarly define the lines $l_B$, $l_C$. Let $A^*=l_B\cap l_C$, $B^*=l_A\cap l_C$, $C^*=l_B\cap l_A$. Consider the point $A^{!}=l_A\cap BC$. Similarly define the points $B^{!}$, $C^{!}$. Let $\pi$ be the pedal circle of $H$ wrt $A^*B^*C^*$. Then the nine~- point circles of the triangles $ABC$, $A^*B^*C^*$, circle $\pi$ and circle $(A^{!}B^{!}C^{!})$ meet at the same point.

\end{theorem}

\begin{theorem}

Let given triangle $ABC$ with orthocenter $H$. Let $AA'$, $BB'$, $CC'$ be altitudes of $ABC$. Let $I_1$, $I_2$ be the incenters of the triangles $C'B'A$, $C'A'B$ respectively. Let $I_1'$, $I_2'$ be $C'$~-- excenter of the triangle $AC'B'$ and respectively $C'$~-- excenter of the triangle $C'BA'$. Then the points $I_1$, $I_2$, $I_1'$, $I_2'$, $A$, $B$ lie on the same conic.

\end{theorem}

\begin{theorem}

Let $AA'$, $BB'$, $CC'$ be altitudes of $ABC$. Let $H$ be the orthocenter of $ABC$ and $A_1B_1C_1$ be the midpoint triangle of $ABC$. Places of the points $A'$, $B'$, $C'$ wrt points $A_1$, $B_1$, $C_1$ are the same as on the picture below. Let the incenters of the triangles $C'C_1H$, $B'B_1H$ and $H$~-- excenter of $HA'A_1$ lie on the circle $\omega_1$. Let $H$~-- excenters of the triangles $C'C_1H$, $B'B_1H$ and the incenter of $HA'A_1$ lie on the circle $\omega_2$. Then

\begin{enumerate}

\item Point $H$ lies on the radical line of $\omega_1$ and $(A'B'C')$.

\item Circles $\omega_1$ and $\omega_2$ are symmetric wrt Euler line of $ABC$.

\item Circles $\omega_1$, $\omega_2$, $(A'B'C')$ are coaxial.

\end{enumerate}

\end{theorem}

\newpage

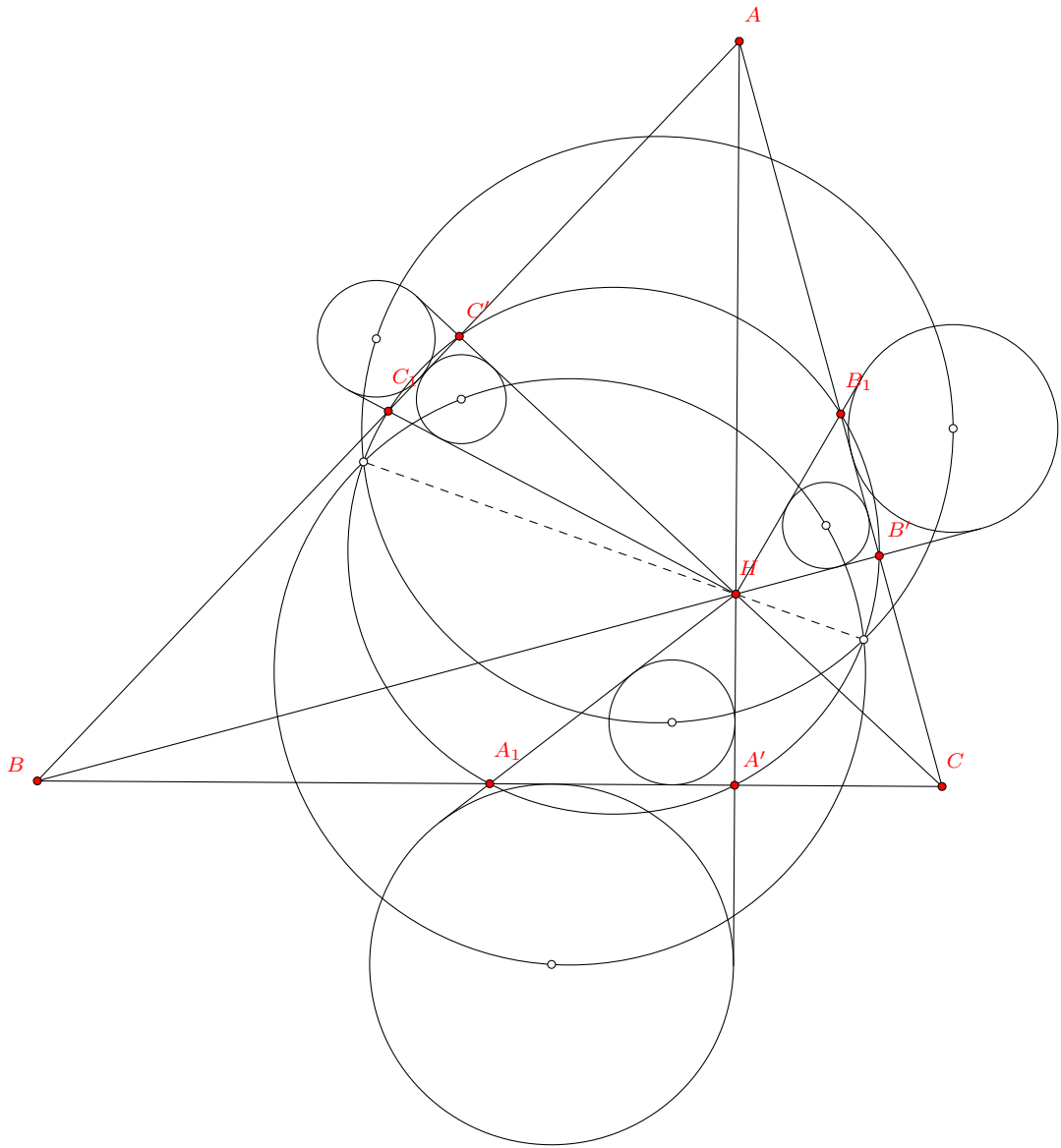
\definecolor{ffffff}{rgb}{1.0,1.0,1.0}
\definecolor{ffqqqq}{rgb}{1.0,0.0,0.0}
\begin{tikzpicture}[line cap=round,line join=round,>=triangle 45,x=1.0cm,y=1.0cm]
%\clip(-4.724556940268088,-10.742057966037688) rectangle (11.327094324715596,12.070899780337687);
\draw (-2.706386413568172,-0.46186748741931893)-- (6.635377852205561,2.0580639526255093);
\draw (-2.706386413568172,-0.46186748741931893)-- (6.681921932919673,9.516054964849946);
\draw (6.681921932919673,9.516054964849946)-- (9.393826605659127,-0.5373829030991867);
\draw (-2.706386413568172,-0.46186748741931893)-- (9.393826605659127,-0.5373829030991867);
\draw(5.586838173417778,4.277909989640748) circle (3.9549167901084714cm);
\draw(4.415531815190315,1.0095242738377244) circle (3.954916790108472cm);
\draw(9.541726491538917,4.292916917068917) circle (1.400798963294004cm);
\draw(2.96381314603573,4.688366445270191) circle (0.6003375308040846cm);
\draw(7.8414551554829375,2.9854838086258595) circle (0.5813039128268213cm);
\draw(5.7819299894785186,0.327807965643392) circle (0.8426332137884553cm);
\draw(1.8264950531208803,5.503139125075294) circle (0.7863004077464303cm);
\draw (9.393826605659127,-0.5373829030991867)-- (2.365316103278257,6.07579976378683);
\draw (1.457599631257715,4.80874406536461)-- (6.635377852205561,2.0580639526255093);
\draw [dash pattern=on 3pt off 3pt] (8.346414575285467,1.4448718856017475)-- (1.6559554133226289,3.842562377876725);
\draw(5.0011849943040465,2.6437171317392365) circle (3.5535603095571036cm);
\draw (6.635377852205561,2.0580639526255093)-- (6.681921932919673,9.516054964849946);
\draw (3.3437200960454776,-0.4996251952592528)-- (6.635377852205561,2.0580639526255093);
\draw (6.635377852205561,2.0580639526255093)-- (6.619288163058765,-0.5200674704303726);
\draw (6.635377852205561,2.0580639526255093)-- (8.0378742692894,4.489336030875379);
\draw (6.635377852205561,2.0580639526255093)-- (8.55409306465493,2.5756354359970137);
\draw(4.171153256948127,-2.9378350886922906) circle (2.432998641022355cm);
\draw (6.604104519134404,-2.9530187326166484)-- (6.619288163058765,-0.5200674704303726);
\draw (2.6783419399926784,-1.016638388927697)-- (3.3437200960454776,-0.4996251952592528);
\draw (8.0378742692894,4.489336030875379)-- (8.32833993889853,4.9928674934934465);
\draw (8.55409306465493,2.5756354359970137)-- (9.906550518632875,2.940459463090973);
\begin{scriptsize}
\draw [fill=ffqqqq] (-2.706386413568172,-0.46186748741931893) circle (1.5pt);
\draw[color=ffqqqq] (-2.9955205499884947,-0.23880705792142481) node {$B$};
\draw [fill=ffqqqq] (6.681921932919673,9.516054964849946) circle (1.5pt);
\draw[color=ffqqqq] (6.862567525933959,9.877346150878514) node {$A$};
\draw [fill=ffqqqq] (9.393826605659127,-0.5373829030991867) circle (1.5pt);
\draw[color=ffqqqq] (9.572251421148247,-0.1871940313459149) node {$C$};
\draw [fill=ffqqqq] (2.937976741989781,5.536978713629455) circle (1.5pt);
\draw[color=ffqqqq] (3.1980426390727326,5.903143104564253) node {$C'$};
\draw [fill=ffqqqq] (6.619288163058765,-0.5200674704303726) circle (1.5pt);
\draw[color=ffqqqq] (6.8883740392217145,-0.16138751805815993) node {$A'$};
\draw [fill=ffqqqq] (8.55409306465493,2.5756354359970137) circle (1.5pt);
\draw[color=ffqqqq] (8.823862535803347,2.9353940764724338) node {$B'$};
\draw [fill=ffqqqq] (6.635377852205561,2.0580639526255093) circle (1.5pt);
\draw[color=ffqqqq] (6.8109544993584485,2.419263810717335) node {$H$};
\draw [fill=ffqqqq] (1.9877677596757506,4.527093738715314) circle (1.5pt);
\draw[color=ffqqqq] (2.217395134138038,4.974108626205075) node {$C_1$};
\draw [fill=ffqqqq] (8.0378742692894,4.489336030875379) circle (1.5pt);
\draw[color=ffqqqq] (8.28192575676049,4.922495599629565) node {$B_1$};
\draw [fill=ffqqqq] (3.3437200960454776,-0.4996251952592528) circle (1.5pt);
\draw[color=ffqqqq] (3.585140338389059,-0.05816146490714019) node {$A_1$};
\draw [fill=ffffff] (1.8264950531208803,5.503139125075294) circle (1.5pt);
\draw [fill=ffffff] (2.96381314603573,4.688366445270191) circle (1.5pt);
\draw [fill=ffffff] (9.541726491538917,4.292916917068917) circle (1.5pt);
\draw [fill=ffffff] (7.8414551554829375,2.9854838086258595) circle (1.5pt);
\draw [fill=ffffff] (4.171153256948127,-2.9378350886922906) circle (1.5pt);
\draw [fill=ffffff] (5.7819299894785186,0.327807965643392) circle (1.5pt);
\draw [fill=ffffff] (8.346414575285467,1.4448718856017475) circle (1.5pt);
\draw [fill=ffffff] (1.6559554133226289,3.842562377876725) circle (1.5pt);
\end{scriptsize}
\end{tikzpicture}
\bigskip
\begin{commentary}

Compare previous theorem with ideas of Section 2.

\end{commentary}
\medskip
\subsection{Analog for cyclic quadrilateral}

\begin{theorem}

Let given convex cyclic quadrilateral $ABCD$, where $AC\perp BD$. Let $AC\cap BD = H$. Let $X$ be the foot of perpendicular from $H$ on $AB$. Let $Y$ be the foot of perpendicular from $H$ on $CD$. Let $M$ and $N$ be the midpoints of segments $AB$ and $CD$ respectively. Let $I_1$, $I_2$, $I_3$, $I_4$ be the incenters of the triangles $HXM$, $ABD$, $HYN$, $CAD$ respectively. Let $W$ be the midpoint of a smaller arc $DA$ of a circle $(ABCD)$. Then $I_1I_2$, $I_3I_4$ and $HW$ are concurrent.

\end{theorem}

\newpage

\definecolor{ffffff}{rgb}{1.0,1.0,1.0}
\definecolor{ffwwqq}{rgb}{1.0,0.4,0.0}
\definecolor{ffqqqq}{rgb}{1.0,0.0,0.0}
\begin{tikzpicture}[line cap=round,line join=round,>=triangle 45,x=1.0cm,y=1.0cm]
%\clip(-0.6177706552032521,-12.490878413173473) rectangle (12.626467071123882,6.33212182810176);
\draw[color=ffqqqq] (3.0831846647205716,1.165670286792168) -- (2.8356040013414914,0.7878741880026369) -- (3.2134001001310226,0.5402935246235565) -- (3.460980763510103,0.9180896234130876) -- cycle; 
\draw (0.9452013275164236,2.566752337703729)-- (8.867626732649128,-2.6250411817596513);
\draw (0.9452013275164236,2.566752337703729)-- (6.1253610446155005,4.983804862054891);
\draw (6.1253610446155005,4.983804862054891)-- (8.867626732649128,-2.6250411817596513);
\draw (8.867626732649128,-2.6250411817596513)-- (1.268565629766738,-2.4274296637177923);
\draw (1.268565629766738,-2.4274296637177923)-- (0.9452013275164236,2.566752337703729);
\draw (1.268565629766738,-2.4274296637177923)-- (6.1253610446155005,4.983804862054891);
\draw (5.068096181207933,-2.526235422738722)-- (2.379460063395984,3.2359746952089066);
\draw (3.535281186065962,3.77527859987931)-- (3.3725586273793966,-2.482143429613157);
\draw [dash pattern=on 3pt off 3pt](3.0274224026585843,3.000310698555249)-- (7.059552659775981,1.517492373625005);
\draw [dash pattern=on 3pt off 3pt](4.03753284020386,-1.850877587024035)-- (7.059552659775981,1.517492373625005);
\draw [dash pattern=on 3pt off 3pt] (3.460980763510103,0.9180896234130876)-- (9.61625495604132,1.9433541966011665);
\draw [shift={(5.1423966037637925,0.3309535537275)}] plot[domain=-0.6707689984678078:1.3625970256191742,variable=\t]({1.0*4.755548800021174*cos(\t r)+-0.0*4.755548800021174*sin(\t r)},{0.0*4.755548800021174*cos(\t r)+1.0*4.755548800021174*sin(\t r)});
\draw(3.0274224026585843,3.000310698555249) circle (0.4875411328398017cm);
\draw(4.998844439377637,2.2753190599398403) circle (1.9781186919636629cm);
\draw(4.03753284020386,-1.8508775870240348) circle (0.648339134921404cm);
\draw(5.515083436985773,-0.20398665225712628) circle (2.333083816703565cm);
\draw (8.124052052367897,3.7543417276475717)-- (8.407440703294913,4.079714451937661);
\draw (9.622449882061154,-0.4031979477648091)-- (10.048252352032058,-0.47297468997066905);
\begin{scriptsize}
\draw [fill=ffqqqq] (8.867626732649128,-2.6250411817596513) circle (1.5pt);
\draw[color=ffqqqq] (9.00665933241711,-2.334123758096226) node {$A$};
\draw [fill=ffqqqq] (1.268565629766738,-2.4274296637177923) circle (1.5pt);
\draw[color=ffqqqq] (1.4263560678311609,-2.1211938911134745) node {$B$};
\draw [fill=ffqqqq] (0.9452013275164236,2.566752337703729) circle (1.5pt);
\draw[color=ffqqqq] (1.0856682806587585,2.8613649962829104) node {$C$};
\draw [fill=ffqqqq] (6.1253610446155005,4.983804862054891) circle (1.5pt);
\draw[color=ffqqqq] (6.302450021736167,4.458338998653547) node {$D$};
\draw [fill=ffqqqq] (3.460980763510103,0.9180896234130876) circle (1.5pt);
\draw[color=ffqqqq] (3.260980763510103,0.9180896234130876) node {$H$};
\draw [fill=ffqqqq] (3.535281186065962,3.77527859987931) circle (1.5pt);
\draw[color=ffqqqq] (3.6834126578483253,4.075065238084593) node {$N$};
\draw [fill=ffqqqq] (5.068096181207933,-2.526235422738722) circle (1.5pt);
\draw[color=ffqqqq] (5.216507700124136,-2.22765882460485) node {$M$};
\draw [fill=ffqqqq] (2.379460063395984,3.2359746952089066) circle (1.5pt);
\draw[color=ffqqqq] (2.5335913761414677,3.542740570627715) node {$Y$};
\draw [fill=ffqqqq] (3.3725586273793966,-2.482143429613157) circle (1.5pt);
\draw[color=ffqqqq] (3.513068764262124,-2.1850728512082997) node {$X$};
\draw [fill=ffwwqq] (5.515083436985773,-0.20398665225712628) circle (1.5pt);
\draw[color=ffwwqq] (5.706246394184463,0.15715568560196627) node {$I_2$};
\draw [fill=ffwwqq] (4.9988444393776374,2.2753190599398394) circle (1.5pt);
\draw[color=ffwwqq] (5.19521471342586,2.6271421426018837) node {$I_4$};
\draw [fill=ffwwqq] (3.0274224026585843,3.000310698555249) circle (1.5pt);
\draw[color=ffwwqq] (3.2149669504862723,3.3723966770415137) node {$I_3$};
\draw [fill=ffwwqq] (4.03753284020386,-1.850877587024035) circle (1.5pt);
\draw[color=ffwwqq] (4.237030312003479,-1.4824042901652201) node {$I_1$};
\draw [fill=ffffff] (7.059552659775981,1.517492373625005) circle (1.5pt);
\draw[color=ffwwqq] (10.01625495604132,1.9433541966011665) node {$W$};
\draw [fill=ffwwqq] (9.61625495604132,1.9433541966011665) circle (1.5pt);
\end{scriptsize}
\end{tikzpicture}
\bigskip
\subsection{Tangent circles}

\begin{theorem}

Let given triangle $ABC$ with orthocenter $H$. Let $AA'$, $BB'$, $CC'$ be altitudes of $ABC$. Let the circle $\omega_1$ is internally tangent to the incircles of the triangles $HA'C$, $HB'A$, $HC'B$. Similarly let the circle $\omega_2$ is internally tangent to the incircles of the triangles $HA'B$, $HB'C$, $HC'A$. Then $H$ lies on the radical line of $\omega_1$ and $\omega_2$.

\end{theorem}

\medskip

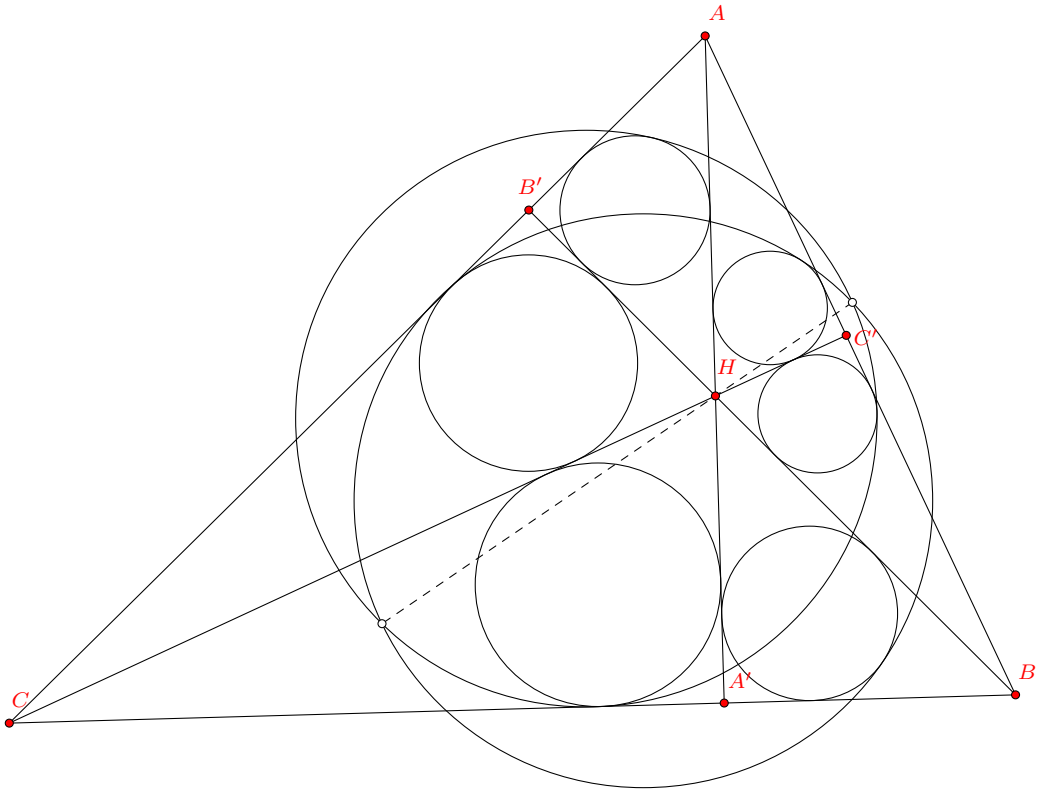
\definecolor{ffffff}{rgb}{1.0,1.0,1.0}
\definecolor{ffqqqq}{rgb}{1.0,0.0,0.0}
\begin{tikzpicture}[line cap=round,line join=round,>=triangle 45,x=1.0cm,y=1.0cm]
%\clip(-2.6663012237062267,-14.61154339413738) rectangle (29.59317758995777,5.114841051427137);
\draw (3.12,-8.08)-- (12.271911605876896,1.0327592590004755);
\draw (12.271911605876896,1.0327592590004755)-- (16.3539933983036,-7.70510230749119);
\draw (3.12,-8.08)-- (16.3539933983036,-7.70510230749119);
\draw (3.12,-8.08)-- (14.126888121929667,-2.9378931753610886);
\draw (12.522515885086076,-7.813641886994741)-- (12.271911605876896,1.0327592590004755);
\draw (16.3539933983036,-7.70510230749119)-- (9.952812224389636,-1.2764188964812164);
\draw(11.459389371920745,-5.131212370658879) circle (3.8045888002400825cm);
\draw(11.34944485498918,-1.2794127333598935) circle (0.9875706728924905cm);
\draw(10.86262420142478,-6.245203585959058) circle (1.6148125441740344cm);
\draw(13.748838932014372,-3.978872224369569) circle (0.7831215009403726cm);
\draw(10.709259725571584,-4.041991493802813) circle (3.823355997556188cm);
\draw(13.126780658661646,-2.5746871947250307) circle (0.752374083303709cm);
\draw(9.948461692996794,-3.305953019387166) circle (1.4351006381370994cm);
\draw(13.645952532605483,-6.624699263554329) circle (1.156653462110199cm);
\draw [dash pattern=on 3pt off 3pt] (14.20827277116278,-2.5009013631383077)-- (8.021770830265389,-6.7614499916760185);
\begin{scriptsize}
\draw [fill=ffqqqq] (3.12,-8.08) circle (1.5pt);
\draw[color=ffqqqq] (3.2692663555520602,-7.77129822844946) node {$C$};
\draw [fill=ffqqqq] (12.271911605876896,1.0327592590004755) circle (1.5pt);
\draw[color=ffqqqq] (12.426368754779528,1.3416735568057363) node {$A$};
\draw [fill=ffqqqq] (16.3539933983036,-7.70510230749119) circle (1.5pt);
\draw[color=ffqqqq] (16.508450547206234,-7.396188009685929) node {$B$};
\draw [fill=ffqqqq] (14.126888121929667,-2.9378931753610886) circle (1.5pt);
\draw[color=ffqqqq] (14.390181076541564,-2.961061305481826) node {$C'$};
\draw [fill=ffqqqq] (12.522515885086076,-7.813641886994741) circle (1.5pt);
\draw[color=ffqqqq] (12.735283052584794,-7.50651454461638) node {$A'$};
\draw [fill=ffqqqq] (9.952812224389636,-1.2764188964812164) circle (1.5pt);
\draw[color=ffqqqq] (9.977119679323508,-0.9531183697476303) node {$B'$};
\draw [fill=ffqqqq] (12.407153139987782,-3.7413047199551643) circle (1.5pt);
\draw[color=ffqqqq] (12.55876059669607,-3.3582368312314475) node {$H$};
\draw [fill=ffffff] (8.021770830265389,-6.7614499916760185) circle (1.5pt);
\draw [fill=ffffff] (14.20827277116278,-2.5009013631383077) circle (1.5pt);
\end{scriptsize}
\end{tikzpicture}

\newpage

\begin{theorem}

Let given triangle $ABC$ with orthocenter $H$. Let $W_A$, $W_B$, $W_C$ be the reflections of $H$ wrt lines $BC$, $AC$, $AB$ respectively. Let $I_{AB}$ be $H$~-- excenter of the triangle $AHW_B$ (and let $(I_{AB})$ be the excircle itself), $I_{CB}$ be $H$~-- excenter of the triangle $CHW_B$, $I_{CA}$ be $H$~-- excenter of the triangle $CHW_A$, $I_{BA}$ be $H$~-- excenter of the triangle $BHW_C$, $I_{BC}$ be $H$~-- excenter of the triangle $BHW_C$, $I_{AC}$ be $H$~-- excenter of the triangle $AHW_C$. Suppose that the circle $\omega_1$ externally touches to $(I_{AB})$, $(I_{BC})$, $(I_{CA})$. Suppose that the circle $\omega_2$ externally touches to $(I_{AC})$, $(I_{BA})$, $(I_{CB})$. Then we have that :

\begin{enumerate}

\item Point $H$ lies on the radical line of $(I_{AB}I_{BC}I_{CA})$ and $(I_{BA}I_{CB}I_{AC})$.

\item Point $H$ lies on the radical line of $\pi_1$ and $\pi_2$.

\end{enumerate}

\end{theorem}
\bigskip
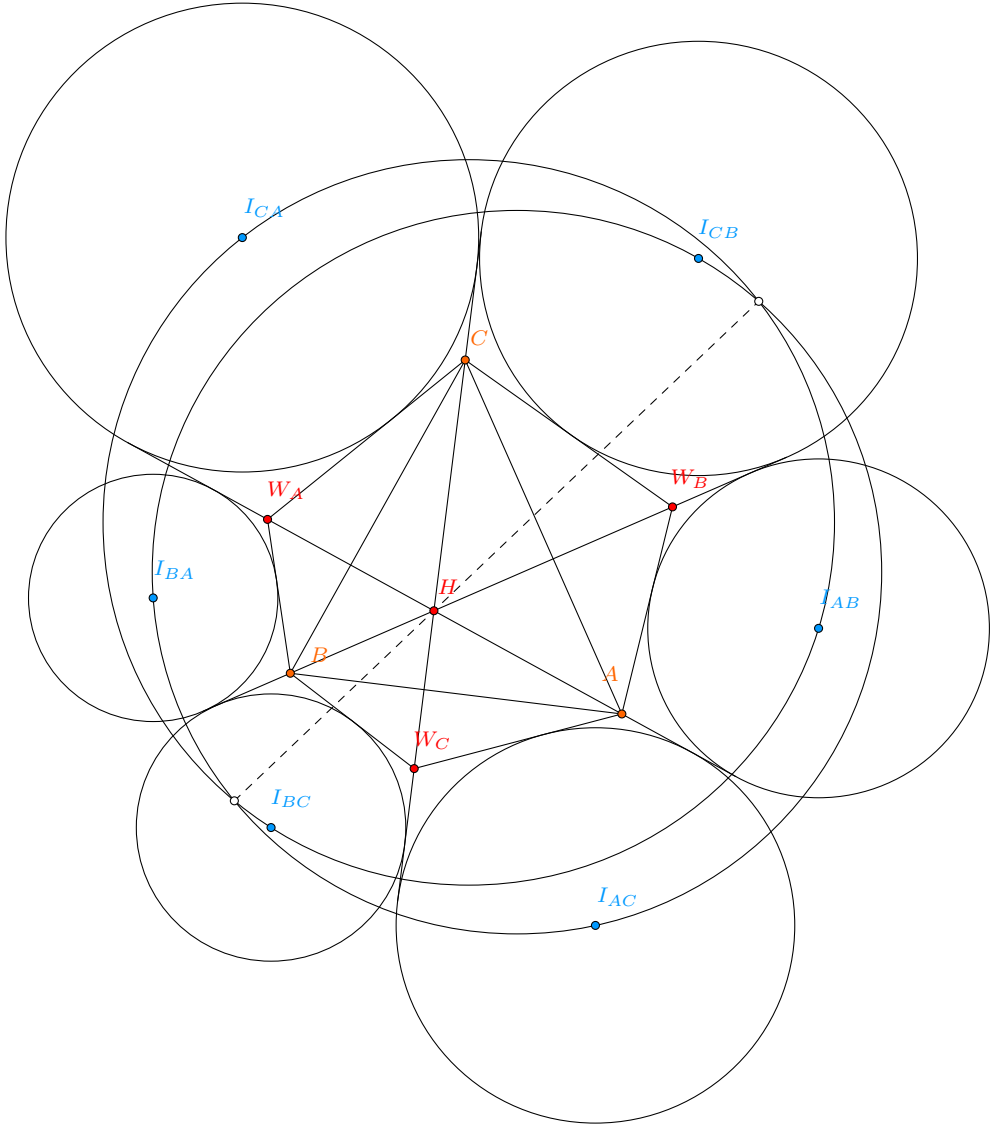
\definecolor{ffffff}{rgb}{1.0,1.0,1.0}
\definecolor{qqccqq}{rgb}{0.0,0.6,1.0}
\definecolor{ffqqqq}{rgb}{1.0,0.0,0.0}
\definecolor{ffwwqq}{rgb}{1.0,0.4,0.0}
\begin{tikzpicture}[line cap=round,line join=round,>=triangle 45,x=1.0cm,y=1.0cm]
%\clip(-7.7900415256925015,-9.232308571287218) rectangle (12.057068346388274,6.378354820983377);
\draw (3.56,-3.7)-- (1.5,0.994);
\draw (-0.8,-3.16)-- (1.5,0.994);
\draw (-0.8,-3.16)-- (3.56,-3.7);
\draw (-0.8,-3.16)-- (0.8289351111548515,-4.424227621046015);
\draw (0.8289351111548515,-4.424227621046015)-- (3.56,-3.7);
\draw (3.56,-3.7)-- (4.226211552294074,-0.9542062637993629);
\draw (4.226211552294074,-0.9542062637993629)-- (1.5,0.994);
\draw (-1.0998017317305564,-1.119946080168445)-- (1.5,0.994);
\draw (-1.0998017317305564,-1.119946080168445)-- (-0.8,-3.16);
\draw(6.147273394120497,-2.564554070752403) circle (2.246598122667491cm);
\draw(4.568310890128029,2.3405531388290886) circle (2.87953453265855cm);
\draw(-1.4314371892359248,2.6170875145023174) circle (3.1087093248602167cm);
\draw(-2.603914035274589,-2.1612031099827744) circle (1.6395236096113004cm);
\draw(-1.0534445552856933,-5.2063301130474455) circle (1.7719747766025078cm);
\draw(3.2127112761492187,-6.504034881461859) circle (2.6213383240084633cm);
\draw(2.180627282238915,-1.8196689901693406) circle (4.796715665326703cm);
\draw(1.5469580307517703,-1.1598499897652341) circle (4.809999479998761cm);
\draw [dash pattern=on 3pt off 3pt] (-1.5362312018491084,-4.851738280811207)-- (5.360441159849584,1.771607228355934);
\draw (5.059042097668666,-4.5299944209528)-- (-2.9372671571109747,-0.10257234921635994);
\draw (0.6112497297168732,-6.181835515619321)-- (1.7106109632199855,2.6944885178502527);
\draw (-1.9450511335722407,-3.6625149840559903)-- (5.725487564265295,-0.29623681670504753);
\begin{scriptsize}
\draw [fill=ffwwqq] (3.56,-3.7) circle (1.5pt);
\draw[color=ffwwqq] (3.404447354290962,-3.172906333497974) node {$A$};
\draw [fill=ffwwqq] (-0.8,-3.16) circle (1.5pt);
\draw[color=ffwwqq] (-0.4211921941437536,-2.91615200138826) node {$B$};
\draw [fill=ffwwqq] (1.5,0.994) circle (1.5pt);
\draw[color=ffwwqq] (1.6841933291558888,1.294619045211045) node {$C$};
\draw [fill=ffqqqq] (0.8289351111548515,-4.424227621046015) circle (1.5pt);
\draw[color=ffqqqq] (1.0679829320925789,-4.045871062671001) node {$W_C$};
\draw [fill=ffqqqq] (4.226211552294074,-0.9542062637993629) circle (1.5pt);
\draw[color=ffqqqq] (4.457140115940784,-0.5796875791898652) node {$W_B$};
\draw [fill=ffqqqq] (-1.0998017317305564,-1.119946080168445) circle (1.5pt);
\draw[color=ffqqqq] (-0.8576745587302648,-0.7337401784556936) node {$W_A$};
\draw [fill=ffqqqq] (1.0881417959804562,-2.3313736472689097) circle (1.5pt);
\draw[color=ffqqqq] (1.273386397780349,-2.0175118390042623) node {$H$};
\draw [fill=qqccqq] (6.147273394120497,-2.564554070752403) circle (1.5pt);
\draw[color=qqccqq] (6.43414847318557,-2.1715644382700905) node {$I_{AB}$};
\draw [fill=qqccqq] (4.568310890128029,2.3405531388290886) circle (1.5pt);
\draw[color=qqccqq] (4.842271614105353,2.7324433050254417) node {$I_{CB}$};
\draw [fill=qqccqq] (-1.4314371892359248,2.6170875145023174) circle (1.5pt);
\draw[color=qqccqq] (-1.1401043240509485,3.014873070346127) node {$I_{CA}$};
\draw [fill=qqccqq] (-2.603914035274589,-2.1612031099827744) circle (1.5pt);
\draw[color=qqccqq] (-2.3211742517556258,-1.7864329401055197) node {$I_{BA}$};
\draw [fill=qqccqq] (-1.0534445552856933,-5.2063301130474455) circle (1.5pt);
\draw[color=qqccqq] (-0.7806482590973511,-4.816134059000142) node {$I_{BC}$};
\draw [fill=qqccqq] (3.2127112761492187,-6.504034881461859) circle (1.5pt);
\draw[color=qqccqq] (3.5071490871348474,-6.125581152759683) node {$I_{AC}$};
\draw [fill=ffffff] (-1.5362312018491084,-4.851738280811207) circle (1.5pt);
\draw [fill=ffffff] (5.360441159849584,1.771607228355934) circle (1.5pt);
\end{scriptsize}
\end{tikzpicture}

\newpage

\definecolor{ffffff}{rgb}{1.0,1.0,1.0}
\definecolor{ffqqqq}{rgb}{1.0,0.0,0.0}
\definecolor{ffwwqq}{rgb}{1.0,0.4,0.0}
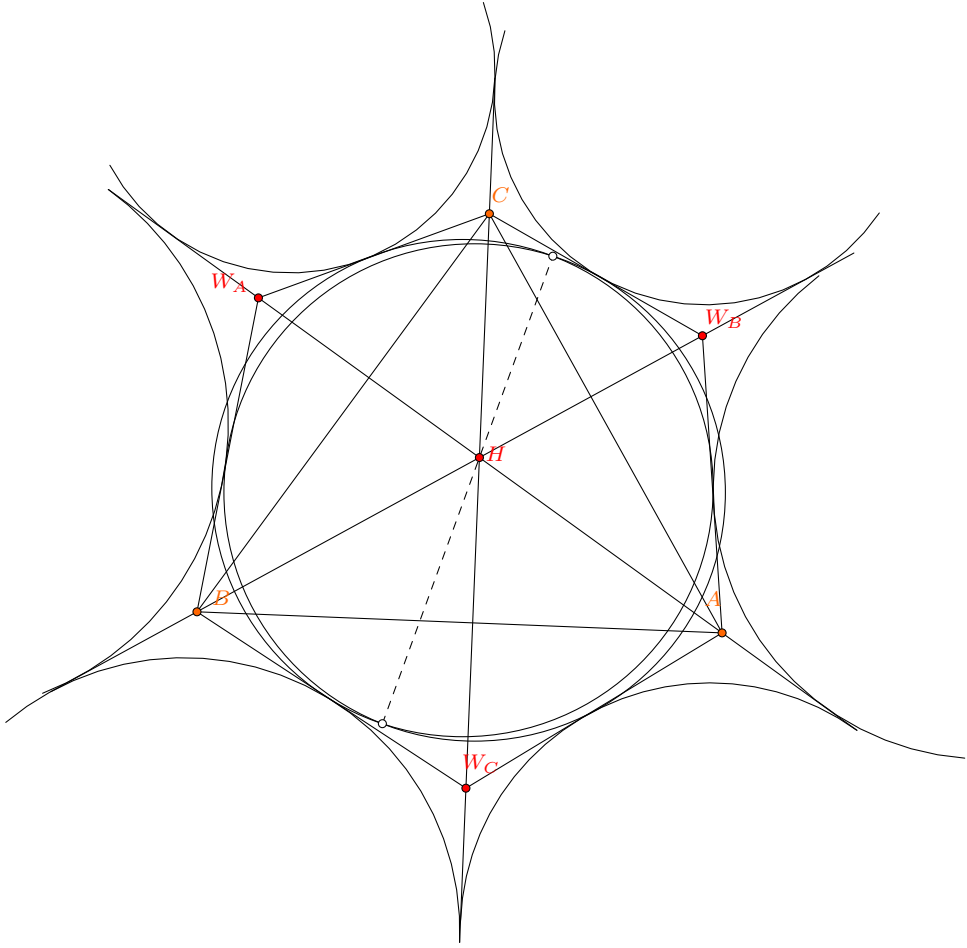
\begin{tikzpicture}[line cap=round,line join=round,>=triangle 45,x=1.0cm,y=1.0cm]
%\clip(-5.000483538965797,-7.815919295336986) rectangle (11.71608293330062,5.332427942875217);
\draw (6.106000061798027,-3.439493936676183)-- (3.0453227052096987,2.118545925594674);
\draw (-0.8,-3.16)-- (3.0453227052096987,2.118545925594674);
\draw (-0.8,-3.16)-- (6.106000061798027,-3.439493936676183);
\draw (-0.8,-3.16)-- (2.7370337380672076,-5.49894893478352);
\draw (2.7370337380672076,-5.49894893478352)-- (6.106000061798027,-3.439493936676183);
\draw (6.106000061798027,-3.439493936676183)-- (5.847498958402081,0.5006159804713985);
\draw (5.847498958402081,0.5006159804713985)-- (3.0453227052096987,2.118545925594674);
\draw (0.008081054900604007,1.0027268457660514)-- (3.0453227052096987,2.118545925594674);
\draw (0.008081054900604007,1.0027268457660514)-- (-0.8,-3.16);
\draw(2.6918916166304383,-1.5196094505620066) circle (3.29685131973788cm);
\draw(2.8517515902461157,-1.577415352074929) circle (3.2975434398112577cm);
\draw [dash pattern=on 3pt off 3pt] (1.6383539545324817,-4.64359520901203)-- (3.88004006687793,1.5557005963490176);
\draw (7.879350896625797,-4.7313470957822235)-- (-1.966142473235625,2.440912032540985);
\draw (2.6542800502119532,-7.543705183304687)-- (3.118026662935404,3.9149839539745197);
\draw (-2.7218734542577208,-4.21832896298651)-- (7.839809524974673,1.5977329478745435);
\draw [shift={(-4.227421199005969,-0.6631872900926004)}] plot[domain=-1.1985815754994675:0.9412099352943185,variable=\t]({1.0*3.84041847724904*cos(\t r)+-0.0*3.84041847724904*sin(\t r)},{0.0*3.84041847724904*cos(\t r)+1.0*3.84041847724904*sin(\t r)});
\draw [shift={(0.42998305236338813,4.023772306246498)}] plot[domain=-2.6531941154819414:0.33953478108870094,variable=\t]({1.0*2.6902441074978842*cos(\t r)+-0.0*2.6902441074978842*sin(\t r)},{0.0*2.6902441074978842*cos(\t r)+1.0*2.6902441074978842*sin(\t r)});
\draw [shift={(5.8797036714335045,3.679370579292234)}] plot[domain=2.8235999407123487:5.688960915183469,variable=\t]({1.0*2.768945836053883*cos(\t r)+-0.0*2.768945836053883*sin(\t r)},{0.0*2.768945836053883*cos(\t r)+1.0*2.768945836053883*sin(\t r)});
\draw [shift={(-0.9714750171416897,-7.396966615663481)}] plot[domain=-0.04044910050696249:2.2730897292136176,variable=\t]({1.0*3.628723193586693*cos(\t r)+-0.0*3.628723193586693*sin(\t r)},{0.0*3.628723193586693*cos(\t r)+1.0*3.628723193586693*sin(\t r)});
\draw [shift={(5.946327076362271,-7.38484481047482)}] plot[domain=0.9412099352943183:3.101143553082832,variable=\t]({1.0*3.282930278209549*cos(\t r)+-0.0*3.282930278209549*sin(\t r)},{0.0*3.282930278209549*cos(\t r)+1.0*3.282930278209549*sin(\t r)});
\draw [shift={(9.564821909251618,-1.5348048787878967)}] plot[domain=2.2282777822463427:4.637675546441167,variable=\t]({1.0*3.576095770722571*cos(\t r)+-0.0*3.576095770722571*sin(\t r)},{0.0*3.576095770722571*cos(\t r)+1.0*3.576095770722571*sin(\t r)});
\begin{scriptsize}
\draw [fill=ffwwqq] (6.106000061798027,-3.439493936676183) circle (1.5pt);
\draw[color=ffwwqq] (5.985306587698288,-2.9934169365848127) node {$A$};
\draw [fill=ffwwqq] (-0.8,-3.16) circle (1.5pt);
\draw[color=ffwwqq] (-0.480739175830376,-2.9717913654693318) node {$B$};
\draw [fill=ffwwqq] (3.0453227052096987,2.118545925594674) circle (1.5pt);
\draw[color=ffwwqq] (3.195607913801306,2.3697247000543755) node {$C$};
\draw [fill=ffqqqq] (2.7370337380672076,-5.49894893478352) circle (1.5pt);
\draw[color=ffqqqq] (2.93610106041554,-5.177599619248352) node {$W_C$};
\draw [fill=ffqqqq] (5.847498958402081,0.5006159804713985) circle (1.5pt);
\draw[color=ffqqqq] (6.136685585506652,0.7261812952778502) node {$W_B$};
\draw [fill=ffqqqq] (0.008081054900604007,1.0027268457660514) circle (1.5pt);
\draw[color=ffqqqq] (-0.3726113202529736,1.2019438598184233) node {$W_A$};
\draw [fill=ffqqqq] (2.914476396640085,-1.1145282004988823) circle (1.5pt);
\draw[color=ffqqqq] (3.1307312004548646,-1.0687411073070392) node {$H$};
\draw [fill=ffffff] (1.6383539545324817,-4.64359520901203) circle (1.5pt);
\draw [fill=ffffff] (3.88004006687793,1.5557005963490176) circle (1.5pt);
\end{scriptsize}
\end{tikzpicture}
\bigskip
\subsection{Conjugations} Note that isogonal conjugation can be considered as a conjugation associated with the incenter.

\begin{theorem}

Let $I$ be the incenter of $ABC$. Let $A'B'C'$ be the cevian triangle of $I$ wrt $ABC$. Let $H_{IAB'}$ be the orthocenter of $IAB'$. Similarly define the points $H_{IA'B}$, $H_{IAC'}$, $H_{IA'C}$, $H_{IBC'}$, $H_{IB'C}$. Let $i_{IAB'}$ be isogonal conjugation of $H_{IAB'}$ wrt $ABC$, similarly define the points $i_{IA'B}$, $i_{IAC'}$, $i_{IA'C}$, $i_{IBC'}$, $i_{IB'C}$. Then $i_{IAB'}i_{IA'B}$, $i_{IAC'}i_{IA'C}$, $i_{IBC'}i_{IB'C}$ are concurrent.

\end{theorem}

\begin{theorem}

Let given triangle $ABC$ with orthocenter $H$. Let $A'$ be the reflection of $H$ wrt $BC$, similarly define the points $B', C'$. Let $A'B'\cap CA = C_B$, $A'B'\cap CB = C_A$. Like the same define the points $A_B$, $A_C$, $B_A$, $B_C$. Let $i_{CC_BB'}(H)$ be isogonal conjugation of $H$ wrt $CC_BB'$. Similarly define the points $i_{AA_BB'}(H)$, $i_{BB_AA'}(H)$, $i_{BB_CC'}(H)$, $i_{CC_AA'}(H)$, $i_{AA_CC'}(H)$. Then we have that :

\begin{enumerate}

\item Lines $i_{AA_BB'}(H)i_{BB_AA'}(H)$, $i_{BB_CC'}(H)i_{CC_BB'}(H)$, $i_{CC_AA'}(H)i_{AA_CC'}(H)$ are concurrent.

\item Let $i_{CC_AC_B}(H)$ be isogonal conjugation of $H$ wrt $CC_AC_B$, similarly define the points $i_{BB_CB_A}(H)$, $i_{AA_BA_C}(H)$. Then the point $H$ is the orthocenter of the triangle $i_{CC_AC_B}(H)i_{BB_CB_A}(H)i_{AA_BA_C}(H)$.

\end{enumerate}

\end{theorem}

\newpage

\section{Combination of different facts}

Here we will build combinations of some well-known constructions from geometry.

\begin{theorem}[\textbf{Pappus and Ceva}]

Consider triangles $A_1B_1C_1$ and $A_2B_2C_2$ and arbitrary points $P$, $Q$. Let $A_1^*B_1^*C_1^*$ be the cevian triangle of $P$ wrt $A_1B_1C_1$ and let $A_2^*B_2^*C_2^*$ be the cevian triangle of $Q$ wrt $A_2B_2C_2$. Let $R=A_1C_2^*\cap C_1^*A_2$, $S=A_1B_2\cap A_2B_1$ and $T = B_1C_2^*\cap C_1^*B_2$. Then by the Pappus's theorem points $R$, $S$, $T$ lie on the same line $l_C$. Similarly define the lines $l_A$, $l_B$. Then the triangle formed by the lines $l_A$, $l_B$, $l_C$ is perspective to the both triangles $A_1B_1C_1$ and $A_2B_2C_2$.

\end{theorem}

\begin{theorem}[\textbf{Two types of tangency}]

Let given triangle $ABC$ and let circle $\omega_A$ is tangent to the lines $AC$, $AB$ and externally tangent to the circle $(ABC)$, similarly define the circles $\omega_B$ and $\omega_C$. Let the circle $\Omega_A$ goes through the points $B$, $C$ and internally tangent to the incircle of $ABC$, similarly define the circles $\Omega_B$ and $\Omega_C$. Let $h_A$ be the external homothety center of the circles $\omega_B$ and $\omega_C$, let $H_A$ be the external homothety center of the circles $\Omega_B$ and $\Omega_C$. Let the line $\mathcal{L}_A$ goes through $h_A$ and $H_A$. Similarly define the lines $\mathcal{L}_B$ and $\mathcal{L}_C$. Let the lines $\mathcal{L}_A$, $\mathcal{L}_B$, $\mathcal{L}_C$ form a triangle $A^{\mathcal{H}}B^{\mathcal{H}}C^{\mathcal{H}}$. Let $\mathcal{H}_A$ be the external homothety center of the circles $\omega_A$ and $\Omega_A$, similarly define $\mathcal{H}_B$, $\mathcal{H}_C$. Then the triangles $A^{\mathcal{H}}B^{\mathcal{H}}C^{\mathcal{H}}$ and $\mathcal{H}_A\mathcal{H}_B\mathcal{H}_C$ are perspective at the incenter of $ABC$.

\end{theorem}

\begin{theorem}[\textbf{Conic and Poncelet point}]

Consider a conic $\mathcal{C}$ and any point $X$. Let given points $A_1$, $B_1$, $C_1$, $D_1$, $A_2$, $B_2$, $C_2$, $D_2$ on $\mathcal{C}$, such that the lines $A_1A_2$, $B_1B_2$, $C_1C_2$, $D_1D_2$ goes through the point $X$. Consider the circle $\Omega_{ABC}$ which goes through the Poncelet points of quadrilaterals $A_1A_2B_1B_2$, $A_1A_2C_1C_2$, $C_1C_2B_1B_2$. Similarly define the circles $\Omega_{BCD}$, $\Omega_{ACD}$, $\Omega_{ABD}$. Then all these circles $\Omega_{ABC}$, $\Omega_{BCD}$, $\Omega_{ACD}$, $\Omega_{ABD}$ goes through the same point.

\end{theorem}

\begin{theorem}[\textbf{Conic and Miquel point}]

Consider a conic $\mathcal{C}$ and any point $X$. Let given points $A_1$, $B_1$, $C_1$, $D_1$, $A_2$, $B_2$, $C_2$, $D_2$ on $\mathcal{C}$, such that the lines $A_1A_2$, $B_1B_2$, $C_1C_2$, $D_1D_2$ goes through the point $X$. Consider the circumcircle $\Omega_{ABC}$ of the triangle $\mathcal{M}(A_1A_2, B_1B_2)\mathcal{M}(A_1A_2, C_1C_2)\mathcal{M}(C_1C_2, B_1B_2)$ (see definition of $\mathcal{M}$ in Section 4). Similarly define a circle $\Omega_{BCD}$. Then we get that $\Omega_{BCD} = \Omega_{ABC}$.

\end{theorem}

\subsection{Combinations with radical lines}

\begin{theorem}[\textbf{Gauss line theorem and radical lines}]

Consider any four lines $l_1$, $l_2$, $l_3$, $l_4$ and let $l$ be the Gauss line for these four lines. Consider any four points $P_1\in l_1$, $P_2\in l_2$, $P_3\in l_3$, $P_4\in l_4$. Let $X_{ij}= l_i\cap l_j$, for any index $i\not= j$. Let $r_1$ be the radical line of circles $(P_1P_2X_{12})$ and $(P_3P_4X_{34})$, let $r_2$ be the radical line of circles $(P_1P_3X_{13})$ and $(P_2P_4X_{24})$, finally let $r_3$ be the radical line of circles $(P_1P_4X_{14})$ and $(P_2P_3X_{23})$. Consider the case when $l_1$, $l_2$, $l_3$ intersect at the same point $P$. Then we get that $P$ lies on $l$.

\end{theorem}

\begin{theorem}[\textbf{Pascal's theorem and radical lines}]

Consider points $A$, $B$, $C$, $D$, $E$, $F$ which lie on the same circle. Let $X = AE\cap BF$, $Y = DB\cap CE$. Let $l_1$ be the radical line of circles $(EFX)$ and $(BCY)$, let $l_2$ be the radical line of circles $(DYE)$ and $(BXA)$, and let $l_3$ be the radical line of circles $(AXF)$ and $(CDY)$. Consider the Pascal line $\mathcal{L}$ of hexagon $AECFBD$ (i.e $\mathcal{L} = XY$). Let the line $\mathcal{W}$ goes through the center of $(ABCDEF)$ and is perpendicular to $\mathcal{L}$. Then we have that :

\begin{enumerate}

\item Lines $\mathcal{W}$, $BE$, $l_3$ are concurrent.

\item Lines $\mathcal{W}$, $l_1$, $l_2$ are concurrent.

\end{enumerate}

\end{theorem}

\begin{theorem}[\textbf{Morley's theorem and radical lines}]

Consider any triangle $ABC$ and it's Morley's triangle $A'B'C'$. Let $l_A$ be the radical line of $(AB'C')$ and $(A'BC)$. Similarly define the lines $l_B$, $l_C$. Then the lines $l_A$, $l_B$, $l_C$ form an equilateral triangle which is perspective to $ABC$.

\end{theorem}

\begin{theorem}[\textbf{Napoleon's theorem and radical lines}]

Consider any triangle $ABC$ and it's Napoleon's triangle $A'B'C'$. Let $l_A$ be the radical line of $(AB'C')$ and $(A'BC)$. Similarly define the lines $l_B$, $l_C$. Then the lines $l_A$, $l_B$, $l_C$ are concurrent.

\end{theorem}

\subsection{On IMO 2011 Problem 6}

Here we will represent some interesting facts which are related to the following IMO 2011 Problem 6 \cite[Problem G8]{IMO}.

\begin{definition}

For any triangle $ABC$ and any point $P$ on $(ABC)$, denote $\omega$ as the circumcircle of triangle formed by the reflections of the tangent line through $P$ to $(ABC)$ wrt sides of $ABC$. By definition let $\otimes(ABC, P)$ be a tangent point of $\omega$ and $(ABC)$.
\end{definition}

\begin{theorem}

For any four lines $l_1$, $l_2$, $l_3$, $l_4$ consider triangles $\triangle_{ijk}$, which are formed by the lines $l_i$, $l_j$, $l_k$. Let $M$ be the Miquel point of the lines $l_1$, $l_2$, $l_3$, $l_4$. Then the points $\otimes(\triangle_{123}, M)$, $\otimes(\triangle_{124}, M)$, $\otimes(\triangle_{234}, M)$, $\otimes(\triangle_{134}, M)$ lie on the circle $\boxtimes(l_1, l_2, l_3, l_4)$.

\end{theorem}

\begin{theorem}

Let given a parabola $\mathcal{K}$ and lines $l_1$, $l_2$, $l_3$, $l_4$, $l_5$ which are tangent to $\mathcal{K}$. Then the circles $\boxtimes(l_1, l_2, l_3, l_4)$, $\boxtimes(l_1, l_2, l_3, l_5)$, $\boxtimes(l_1, l_2, l_5, l_4)$ goes through the same point.

\end{theorem}
\bigskip
\section{Blow-ups}

\begin{theorem}[\textbf{Blow-up of Desargues's theorem}]

Consider circles $\omega_A$, $\omega_B$, $\omega_C$, $\omega_{A'}$, $\omega_{B'}$, $\omega_{C'}$ with centers at points $A$, $B$, $C$, $A'$, $B'$, $C'$ respectively. Let $H_{AB}$ be the external homothety center of the circles $\omega_A$, $\omega_B$. Similarly define the points $H_{BC}$, $H_{CA}$, $H_{A'B'}$, $H_{B'C'}$, $H_{C'A'}$, $H_{AB'}$, $H_{AA'}$, $H_{BB'}$, $H_{CC'}$. Let the triangle $\triangle^H$ formed by the lines $H_{AB}H_{A'B'}$, $H_{BC}H_{B'C'}$, $H_{CA}H_{C'A'}$. Then we have that :

\begin{enumerate}

\item Triangle $\triangle^H$ is perspective to $H_{AA'}H_{BB'}H_{CC'}$, to $ABC$ and to $A'B'C'$.

\item All these perspective centers lie on the same line.

\end{enumerate}

\end{theorem}

\begin{theorem}

Consider points $A$, $B$, $C$, $D$ and let $P = AB\cap CD$ and $Q= BC\cap AD$. Consider segments $A_1A_2$, $B_1B_2$, $C_1C_2$ and $D_1D_2$. Let given that the segments $A_1A_2$, $B_1B_2$, $C_1C_2$, $D_1D_2$ are parallel to each other and that the points $A$, $B$, $C$, $D$ are the midpoints of the segments $A_1A_2$, $B_1B_2$, $C_1C_2$, $D_1D_2$ respectively. Let also given that $P = A_1B_1\cap C_1D_1 = A_2B_2\cap C_2D_2$ and $Q= B_1C_1\cap A_1D_1 = B_2C_2\cap A_2D_2$. Let the circle $\Omega_{PBC}$ goes through the points $\mathcal{M}(P, B_1B_2)$, $\mathcal{M}(P, C_1C_2)$, $\mathcal{M}(B_1B_2, C_1C_2)$ and let the line $\mathcal{L}_{PBC}$ goes through the center of $\Omega_{PBC}$ and is perpendicular to $DA$ (see definition of $\mathcal{M}$ in Section 4). Similarly define the lines $\mathcal{L}_{PDA}$, $\mathcal{L}_{QAB}$, $\mathcal{L}_{QCD}$. Consider the case when $\mathcal{L}_{PBC}$, $\mathcal{L}_{PDA}$, $\mathcal{L}_{QAB}$ intersect at the same point $W$. Then $\mathcal{L}_{QCD}$ also goes through $W$.

\end{theorem}

\begin{theorem}[\textbf{Blow-up of Miquel point}]

Consider any four points $A$, $B$, $C$, $D$. Let $P= AB\cap CD$, $Q= BC\cap DA$. Consider circles $\omega_A$, $\omega_B$, $\omega_C$, $\omega_D$ with centers at points $A$, $B$, $C$, $D$ and radii $r_A$, $r_B$, $r_C$, $r_D$ respectively. Let given that $\frac{|PA|}{|PB|} = \frac{r_A}{r_B}$, $\frac{|PC|}{|PD|} = \frac{r_C}{r_D}$, $\frac{|QB|}{|QC|} = \frac{r_B}{r_C}$, $\frac{|QD|}{|QA|} = \frac{r_D}{r_A}$. Let the circle $\Omega_{PBC}$ goes through $P$ and internally touches the circles $\omega_B$, $\omega_C$, let the circle $\Omega_{QAB}$ goes through $Q$ and internally touches the circles $\omega_A$, $\omega_B$. Similarly define the circles $\Omega_{QCD}$ and $\Omega_{PDA}$. Then there exists a circle $\omega_M$ which is tangent to $\Omega_{PBC}$, $\Omega_{QAB}$, $\Omega_{QCD}$ and $\Omega_{PDA}$.

\end{theorem}

\begin{theorem}[\textbf{Blow-up of Pascal's theorem}]

Consider two circles $\Omega_1$ and $\Omega_2$. Let circles $\omega_A$, $\omega_B$, $\omega_C$, $\omega_D$, $\omega_E$, $\omega_F$ are tangent to both $\Omega_1$ and $\Omega_2$. Consider two external tangents $\mathcal{L}_{AB}^{(1)}$, $\mathcal{L}_{AB}^{(2)}$ to the circles $\omega_A$, $\omega_B$, similarly consider two external tangents $\mathcal{L}_{CD}^{(1)}$, $\mathcal{L}_{CD}^{(2)}$ to the circles $\omega_C$, $\omega_D$. Let the lines $\mathcal{L}_{AB}^{(1)}$, $\mathcal{L}_{AB}^{(2)}$, $\mathcal{L}_{CD}^{(1)}$, $\mathcal{L}_{CD}^{(2)}$ form a convex quadrilateral $\Box_{AB, CD}$ and let $P_{AB, CD}$ be the intersection of the diagonals of $\Box_{AB, CD}$. Similarly define $\Box_{BC, DE}$, $\Box_{CD, EF}$ and the points $P_{BC, DE}$, $P_{CD, EF}$. Then $P_{AB, CD}$, $P_{BC, DE}$, $P_{CD, EF}$ are collinear.

\end{theorem}

\begin{theorem}[\textbf{Blow-up of Brianchon's theorem}]

Consider two circles $\Omega_1$ and $\Omega_2$. Consider two circles $\omega_A$, $\omega_A'$ which tangent to both $\Omega_1$ and $\Omega_2$ and also tangent to each other. Similarly define pairs of circles $\omega_B$, $\omega_B'$; $\omega_C$, $\omega_C'$; $\omega_D$, $\omega_D'$; $\omega_E$, $\omega_E'$; $\omega_F$, $\omega_F'$. Consider two external tangents $\mathcal{L}_{A}^{(1)}$, $\mathcal{L}_{A}^{(2)}$ to the circles $\omega_A$, $\omega_A'$, similarly consider two external tangents $\mathcal{L}_{B}^{(1)}$, $\mathcal{L}_{B}^{(2)}$ to the circles $\omega_B$, $\omega_B'$. Let the lines $\mathcal{L}_{A}^{(1)}$, $\mathcal{L}_{A}^{(2)}$, $\mathcal{L}_{B}^{(1)}$, $\mathcal{L}_{B}^{(2)}$ form a convex quadrilateral $\Box_{A,B}$ and let $P_{A, B}$ be the intersection of the diagonals of $\Box_{A, B}$. Similarly define $\Box_{B, C}$, $\Box_{C, D}$, $\Box_{D, E}$, $\Box_{E, F}$, $\Box_{F, A}$ and $P_{B, C}$, $P_{C, D}$, $P_{D, E}$, $P_{E, F}$, $P_{F, A}$. Then $P_{A, B}P_{D, E}$, $P_{B, C}P_{E, F}$, $P_{C, D}P_{F, A}$ are concurrent.

\end{theorem}
\bigskip
\section{Generalizations of Feuerbach's theorem}

\subsection{The Feuerbach point for a set of coaxial circles} Consider any triangle $ABC$ and let given three circles $\omega_A$, $\omega_B$, $\omega_C$, such that $\omega_A$ is tangent to $BC$, $\omega_B$ is tangent to $AC$, $\omega_C$ is tangent to $AB$ and also given that the circles $\omega_A$, $\omega_B$, $\omega_C$, $(ABC)$ are coaxial. Then the set of circles $\omega_A$, $\omega_B$, $\omega_C$ can be seen as a generalized incircle for $ABC$.

\bigskip
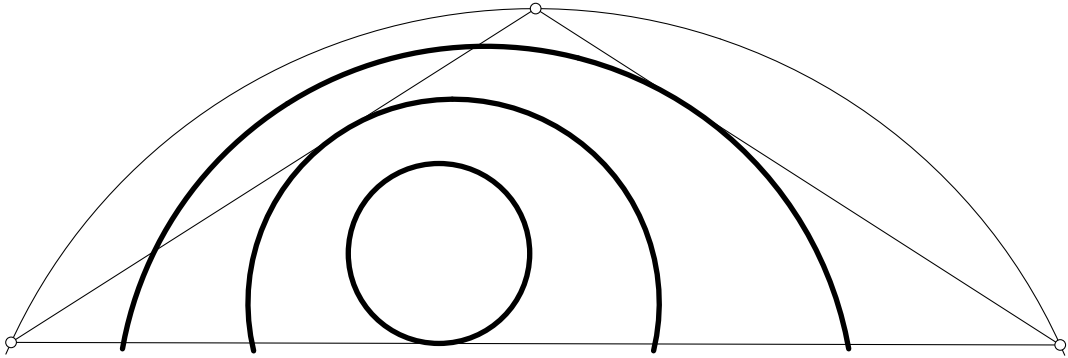
\definecolor{ffffff}{rgb}{1.0,1.0,1.0}
\begin{tikzpicture}[line cap=round,line join=round,>=triangle 45,x=1.0cm,y=1.0cm]
%\clip(-13.63541644805478,-7.155517424305198) rectangle (2.79537199509036,2.8469869152542153);
\draw [line width=2.0pt] (-7.383273107066006,-3.8712030420752663) circle (1.193043622741855cm);
\draw (-13.011318820725373,-5.050325953521931)-- (-6.112203231338463,-0.6237341581992785);
\draw (-6.112203231338463,-0.6237341581992785)-- (0.7869123580484471,-5.084464193588687);
\draw (-13.011318820725373,-5.050325953521931)-- (0.7869123580484471,-5.084464193588687);
\draw [shift={(-6.770342547431304,-5.972058435324329)},line width=2.0pt]  plot[domain=2.138381011587314:2.9682900320739085,variable=\t]({1.0*4.848646581051996*cos(\t r)+-0.0*4.848646581051996*sin(\t r)},{0.0*4.848646581051996*cos(\t r)+1.0*4.848646581051996*sin(\t r)});
\draw [shift={(-6.770342547431294,-5.972058435324332)},line width=2.0pt]  plot[domain=1.3169878454064106:2.138381011587315,variable=\t]({1.0*4.8486465810520025*cos(\t r)+-0.0*4.8486465810520025*sin(\t r)},{0.0*4.8486465810520025*cos(\t r)+1.0*4.8486465810520025*sin(\t r)});
\draw [shift={(-6.770342547431295,-5.972058435324336)},line width=2.0pt]  plot[domain=0.17337808514020336:1.3169878454064106,variable=\t]({1.0*4.848646581052009*cos(\t r)+-0.0*4.848646581052009*sin(\t r)},{0.0*4.848646581052009*cos(\t r)+1.0*4.848646581052009*sin(\t r)});
\draw [shift={(-7.1911135633897265,-4.529841133047875)},line width=2.0pt]  plot[domain=1.5772316122459031:3.378141643753832,variable=\t]({1.0*2.704953898322903*cos(\t r)+-0.0*2.704953898322903*sin(\t r)},{0.0*2.704953898322903*cos(\t r)+1.0*2.704953898322903*sin(\t r)});
\draw [shift={(-7.191113563389729,-4.529841133047876)},line width=2.0pt]  plot[domain=-0.23587762648541233:1.5772316122459022,variable=\t]({1.0*2.704953898322901*cos(\t r)+-0.0*2.704953898322901*sin(\t r)},{0.0*2.704953898322901*cos(\t r)+1.0*2.704953898322901*sin(\t r)});
\draw [shift={(-6.1199568179709125,-8.20129361694651)}] plot[domain=2.009004644801402:2.7351969514000345,variable=\t]({1.0*7.577563425597659*cos(\t r)+-0.0*7.577563425597659*sin(\t r)},{0.0*7.577563425597659*cos(\t r)+1.0*7.577563425597659*sin(\t r)});
\draw [shift={(-6.119956817970907,-8.20129361694654)}] plot[domain=1.2436690536293287:2.009004644801401,variable=\t]({1.0*7.577563425597687*cos(\t r)+-0.0*7.577563425597687*sin(\t r)},{0.0*7.577563425597687*cos(\t r)+1.0*7.577563425597687*sin(\t r)});
\draw [shift={(-6.119956817970911,-8.201293616946517)}] plot[domain=0.40422690663936983:1.2436690536293273,variable=\t]({1.0*7.577563425597666*cos(\t r)+-0.0*7.577563425597666*sin(\t r)},{0.0*7.577563425597666*cos(\t r)+1.0*7.577563425597666*sin(\t r)});
\begin{scriptsize}
\draw [fill=ffffff] (-13.011318820725373,-5.050325953521931) circle (2pt);
\draw [fill=ffffff] (-6.112203231338463,-0.6237341581992785) circle (2pt);
\draw [fill=ffffff] (0.7869123580484471,-5.084464193588687) circle (2pt);
\end{scriptsize}
\end{tikzpicture}
\bigskip
\begin{theorem}

Consider any triangle $ABC$ and let given three circles $\omega_A$, $\omega_B$, $\omega_C$ (with centers at $O_A$, $O_B$, $O_C$ respectively), such that $\omega_A$ is tangent to $BC$ at $A'$, $\omega_B$ is tangent to $AC$ at $B'$, $\omega_C$ is tangent to $AB$ at $C'$ and also given that the circles $\omega_A$, $\omega_B$, $\omega_C$, $(ABC)$ are coaxial. Then we have that :
\begin{enumerate}

\item\textbf{(Gergonne point)} Lines $AA'$, $BB'$, $CC'$ are concurrent.

\item\textbf{(Feuerbach point)} Pedal circles of $O_A$, $O_B$, $O_C$ wrt $ABC$, the nine~- point circle of $ABC$ and the circle $(A'B'C')$ goes through the same point.

\end{enumerate}

\end{theorem}

\subsection{The Feuerbach point for a set of three conics} We can look on the set of three confocal conics (same situation as \cite[Problem 11.18]{Ak}) as on the generalized circle. Also, we can look at a construction with three conics and their three tangents (see picture below), as on a generalized incircle construction.

\newpage
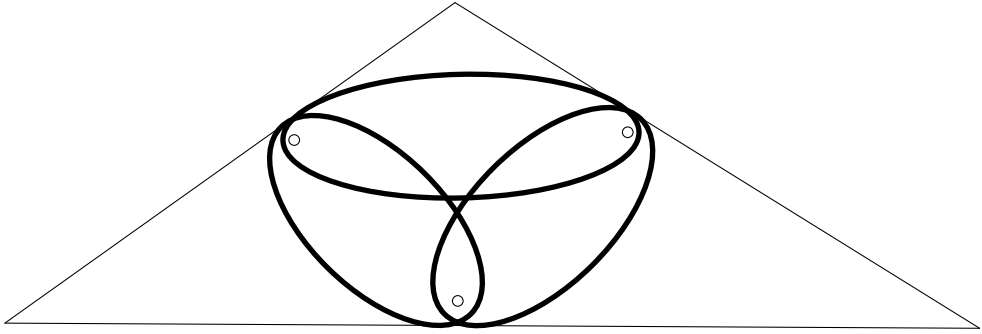
\definecolor{ffwwqq}{rgb}{0.0,0.0,0.0}
\definecolor{ffqqqq}{rgb}{0.0,0.0,0.0}
\definecolor{qqqqcc}{rgb}{0.0,0.0,0.0}
\definecolor{ffffff}{rgb}{1.0,1.0,1.0}
\begin{tikzpicture}[line cap=round,line join=round,>=triangle 45,x=1.0cm,y=1.0cm]
%\clip(-5.253644230726421,-4.431715408530769) rectangle (12.560786688418593,6.491552557433882);
\draw (-0.587235394223108,-1.0126292477674714)-- (5.332800372863854,3.2375557879493955);
\draw (5.332800372863854,3.2375557879493955)-- (12.232241855893378,-1.0809455009203504);
\draw (-0.587235394223108,-1.0126292477674714)-- (12.232241855893378,-1.0809455009203504);
\draw [color=ffqqqq][rotate around={1.3476726290698011:(5.4111686153144705,1.4668181334445674)},line width=2.0pt] (5.4111686153144705,1.4668181334445674) ellipse (2.3419699480147997cm and 0.8198164529236314cm);
\draw [color=ffwwqq][rotate around={-44.758575499605854:(4.293644300695104,0.3485878799609866)},line width=2.0pt] (4.293644300695104,0.3485878799609866) ellipse (1.7588691556845897cm and 0.8937928584691871cm);
\draw [color=qqqqcc][rotate around={45.01809112536193:(6.486829668874906,0.40018402053548174)},line width=2.0pt] (6.486829668874906,0.40018402053548174) ellipse (1.8274206741832264cm and 0.9166061459868521cm);
\begin{scriptsize}
\draw [fill=ffffff] (5.369305354255545,-0.7180462329481028) circle (2pt);
\draw [fill=ffffff] (3.2179832471346654,1.415221992870073) circle (2pt);
\draw [fill=ffffff] (7.6043539834942635,1.5184142740190605) circle (2pt);
\end{scriptsize}
\end{tikzpicture}

\medskip
\begin{theorem}

Consider points $A$, $B$, $C$. Let given set of conics $\mathcal{K}_{AB}$, $\mathcal{K}_{BC}$, $\mathcal{K}_{AC}$, where conic $\mathcal{K}_{AB}$ has foci $A$ and $B$, similarly conic $\mathcal{K}_{BC}$ has foci $B$ and $C$ and conic $\mathcal{K}_{CA}$ has foci $C$ and $A$. Consider the tangent line $l_{B}$ to conics $\mathcal{K}_{AB}$, $\mathcal{K}_{BC}$, similarly define the lines $l_A$, $l_C$. Let the lines $l_A$, $l_B$, $l_C$ form a triangle $A'B'C'$. Let $K_A$ be the foot of perpendicular from $A$ to $l_A$. Similarly define $K_B$ and $K_C$. Then

\begin{enumerate}

\item\textbf{(Gergonne point)} Triangles $A'B'C'$ and $K_AK_BK_C$ are perspective.

\item\textbf{(Another incenter)} Lines $A'A$, $B'B$, $C'C$ are concurrent at $I$.

\item\textbf{(Feuerbach point)} Pedal circle of point $I$ wrt $A'B'C'$, circle $(K_AK_BK_C)$, the nine~- point circle of $A'B'C'$ and the cevian circle of $I$ wrt $ABC$ goes through the same point.

\end{enumerate}

\end{theorem}
\bigskip
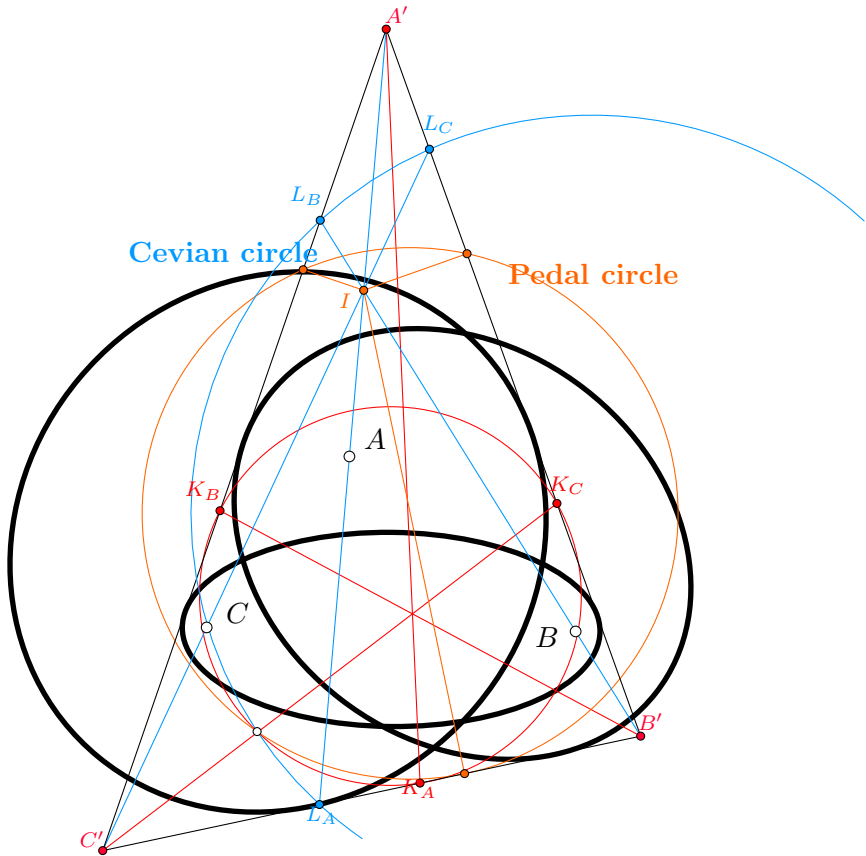
\definecolor{qqccqq}{rgb}{1.0,0.4,0.0}
\definecolor{ffwwqq}{rgb}{1.0,0.4,0.0}
\definecolor{ffqqqq}{rgb}{1.0,0.0,0.0}
\definecolor{wwqqcc}{rgb}{0.0,0.6,1.0}
\definecolor{ffqqtt}{rgb}{1.0,0.0,0.2}
\definecolor{qqqqcc}{rgb}{0.0,0.0,0.0}
\definecolor{ffffff}{rgb}{1.0,1.0,1.0}
\begin{tikzpicture}[line cap=round,line join=round,>=triangle 45,x=1.0cm,y=1.0cm]
%\clip(-9.214920709033631,-3.4893483410637423) rectangle (14.415699188531644,9.182553880818038);
\draw [rotate around={50.37099469553565:(-1.2919157327661936,1.8485485389519756)},line width=2.0pt] (-1.2919157327661936,1.8485485389519756) ellipse (3.7030695692984836cm and 3.398096009400589cm);
\draw [rotate around={-37.945748266635086:(1.134329192017107,1.8219838134981436)},line width=2.0pt] (1.134329192017107,1.8219838134981436) ellipse (3.2197614268874837cm and 2.6092628907875492cm);
\draw [rotate around={-0.6273009872255816:(0.1957088926483848,0.6885555274679991)},line width=2.0pt] (0.1957088926483848,0.6885555274679991) ellipse (2.7467075349448997cm and 1.2872576132413915cm);
\draw [color=ffqqqq] (-2.0551159617903325,2.262782244673629)-- (3.4773406929018877,-0.7264960514806054);
\draw [color=ffqqqq] (-3.5977508857855787,-2.244370274460297)-- (2.3757419973622342,2.3594409050073386);
\draw [color=ffqqqq] (0.1305957478720172,8.648827241497216)-- (0.5758569739244453,-1.3489738380152265);
\draw [color=ffqqqq] (0.18605234050955832,1.1312143720513117) circle (2.510633627111844cm);
%\draw [color=qqqqcc] (0.5115790436095332,1.3611847071662013) circle (2.9034755727223933cm);
\draw [color=qqccqq] (0.4452414350365527,2.227251302272776) circle (3.5241403855928812cm);
\draw [color=wwqqcc](-0.7369159484922225,6.114198038678601)-- (3.4773406929018877,-0.7264960514806054);
\draw [color=wwqqcc](-3.5977508857855787,-2.244370274460297)-- (0.6990223055370557,7.0564792036407855);
\draw [color=wwqqcc](-0.7473547696385957,-1.6328527072511434)-- (0.1305957478720172,8.648827241497216);
\draw (0.1305957478720172,8.648827241497216)-- (-3.5977508857855787,-2.244370274460297);
\draw (-3.5977508857855787,-2.244370274460297)-- (3.4773406929018877,-0.7264960514806054);
\draw (0.1305957478720172,8.648827241497216)-- (3.4773406929018877,-0.7264960514806054);
\draw [color=qqccqq] (-0.16509420707808867,5.186001694703164)-- (-0.961329361237006,5.4585241246195855);
\draw [color=qqccqq] (-0.16509420707808867,5.186001694703164)-- (1.1935955502823705,5.671018364018009);
\draw [color=qqccqq] (-0.16509420707808867,5.186001694703164)-- (1.1611115325210521,-1.2234146422797065);
\draw [shift={(2.835242085620546,2.2358523018691048)},color=wwqqcc]  plot[domain=1.9879320949142947:3.9653690328884137,variable=\t]({1.0*5.27274868307838*cos(\t r)+-0.0*5.27274868307838*sin(\t r)},{0.0*5.27274868307838*cos(\t r)+1.0*5.27274868307838*sin(\t r)});
\draw [color=qqccqq](1.6022114669638943,5.671233488609455) node[anchor=north west] {$\textbf{Pedal circle}$};
%\draw [color=qqqqcc](2.2110291156129698,3.391706943667593) node[anchor=north west] {$\textbf{Nine-point circle}$};
\draw [color=wwqqcc](-3.395756671945957,5.968563037949698) node[anchor=north west] {$\textbf{Cevian circle}$};
\draw [shift={(2.835242085620548,2.2358523018691128)},color=wwqqcc]  plot[domain=0.8240818913943492:1.9879320949142956,variable=\t]({1.0*5.272748683078374*cos(\t r)+-0.0*5.272748683078374*sin(\t r)},{0.0*5.272748683078374*cos(\t r)+1.0*5.272748683078374*sin(\t r)});
\draw [shift={(2.8352420856205702,2.2358523018691345)},color=wwqqcc]  plot[domain=3.965369032888414:4.102560690430101,variable=\t]({1.0*5.272748683078419*cos(\t r)+-0.0*5.272748683078419*sin(\t r)},{0.0*5.272748683078419*cos(\t r)+1.0*5.272748683078419*sin(\t r)});
\draw (-0.29503422882624764,3.476658243479091) node[anchor=north west] {$A$};
\draw (1.953504765518719,0.8573264992912372) node[anchor=north west] {$B$};
\draw (-2.107328624804891,1.168814598600063) node[anchor=north west] {$C$};
\begin{scriptsize}
\draw [fill=ffffff] (-2.2305360321349164,0.7151202529218306) circle (2pt);
\draw [fill=ffffff] (-0.35329543339747155,2.9819768249821204) circle (2pt);
\draw [fill=ffffff] (2.621953817431686,0.6619908020141675) circle (2pt);
\draw [fill=ffqqqq] (0.5758569739244453,-1.3489738380152265) circle (1.5pt);
\draw[color=ffqqqq] (0.5544787692887413,-1.422200045650625) node {$K_A$};
\draw [fill=ffqqqq] (2.3757419973622342,2.3594409050073386) circle (1.5pt);
\draw[color=ffqqqq] (2.508358664953216,2.5988281454269453) node {$K_C$};
\draw [fill=ffqqqq] (-2.0551159617903325,2.262782244673629) circle (1.5pt);
\draw[color=ffqqqq] (-2.2772312244278883,2.5280353955840305) node {$K_B$};
\draw [fill=ffqqqq] (0.1305957478720172,8.648827241497216) circle (1.5pt);
\draw[color=ffqqtt] (0.27130776991707833,8.842748681572047) node {$A'$};
\draw [fill=ffqqtt] (-3.5977508857855787,-2.244370274460297) circle (1.5pt);
\draw[color=ffqqtt] (-3.735561871191953,-2.0876518941740256) node {$C'$};
\draw [fill=ffqqtt] (3.4773406929018877,-0.7264960514806054) circle (1.5pt);
\draw[color=ffqqtt] (3.6127255625027015,-0.5302113976298962) node {$B'$};
\draw [fill=wwqqcc] (-0.7369159484922225,6.114198038678601) circle (1.5pt);
\draw[color=wwqqcc] (-0.9180104274439062,6.435795186912937) node {$L_B$};
\draw [fill=wwqqcc] (0.6990223055370557,7.0564792036407855) circle (1.5pt);
\draw[color=wwqqcc] (0.8234912186918211,7.384418034807998) node {$L_C$};
\draw [fill=wwqqcc] (-0.7473547696385957,-1.6328527072511434) circle (1.5pt);
\draw[color=wwqqcc] (-0.7197907278837421,-1.7761637948651998) node {$L_A$};
\draw [fill=ffwwqq] (-0.16509420707808867,5.186001694703164) circle (1.5pt);
\draw[color=ffwwqq] (-0.40830262857491284,5.0482572899918035) node {$I$};
\draw [fill=qqccqq] (-0.961329361237006,5.4585241246195855) circle (1.5pt);
\draw [fill=qqccqq] (1.1935955502823705,5.671018364018009) circle (1.5pt);
\draw [fill=qqccqq] (1.1611115325210521,-1.2234146422797065) circle (1.5pt);
\draw [fill=ffffff] (-1.5673882298449362,-0.6656516151106914) circle (1.5pt);
\end{scriptsize}
\end{tikzpicture}

\section{Isogonal conjugacy for isogonal lines}

Consider any situation with a triangle $ABC$ and points $P, Q$, such that $\angle BAP = \angle QAC$. Next we will see that it's natural to consider isogonal conjugations of $P$, $Q$ wrt $ABC$.

\begin{theorem}

Let the incircle of $ABC$ is tangent to the sides of $ABC$ at $A_1B_1C_1$. Let $A$~-- excircle of $ABC$ is tangent to the sides of $ABC$ at $A_2B_2C_2$. Let $A_1A_1^*$ be $A_1$~-- altitude of $A_1B_1C_1$. Similarly let $A_2A_2^*$ be $A_2$~-- altitude of $A_2B_2C_2$. Consider isogonal conjugations $B'$, $C'$ of $B$, $C$ wrt $AA_1^*A_2^*$ respectively. Then the midpoint of $BC$ lies on $B'C'$.

\end{theorem}

\begin{theorem}

Consider any four points $A$, $B$, $C$, $D$. Let $P= AB\cap CD$, $Q= BC\cap DA$. Consider a point $N$, such that $\angle DNQ = \angle PNB$. Consider isogonal conjugations $C_1$, $C_2$ of $C$ wrt triangles $NBD$ and $NPQ$ respectively. Consider isogonal conjugations $Q_1$, $Q_2$ of $Q$ wrt triangles $NAC$ and $NBD$ respectively. Consider isogonal conjugations $D_1$, $D_2$ of $D$ wrt triangles $NAC$ and $NPQ$ respectively. Then all the points $C_1$, $C_2$, $Q_1$, $Q_2$, $D_1$, $D_2$ lie on the same conic.

\end{theorem}

\begin{theorem}

Consider any triangle $ABC$ and any two points $P$, $Q$ on it's circumcircle. Consider a point $A'$ on $BC$, such that $BC$ is the external angle bisector of $\angle PA'Q$. Similarly define $B'$ on $AC$ and $C'$ on $AB$. Consider isogonal conjugations $B^A$, $C^A$ of $B$, $C$ wrt $PA'Q$ respectively. Similarly define the points $A^B$, $A^C$, $B^C$, $C^B$. Then all the points $B^A$, $C^A$, $A^B$, $A^C$, $B^C$, $C^B$ lie on the same conic.

\end{theorem}

\begin{commentary}

Note that the previous three theorems are related to \cite[Problem 4.5.7]{Ak} to \cite[Problem 4.12.3]{Ak} and to \cite[Problem 3.13]{Ak}.

\end{commentary}
\bigskip
\section{Fun with some lines}

In this section we will construct some nice theorems which includes some famous lines.

\begin{theorem}[\textbf{Gauss lines}]

Consider any four lines $l_1$, $l_2$, $l_3$, $l_4$. Let $l$ be the Gauss line of the complete quadrilateral formed by these four lines. Let $g_1$ be the Gauss line of the complete quadrilateral formed by four lines $l$, $l_2$, $l_3$, $l_4$. Similarly define the lines $g_2$, $g_3$, $g_4$. Then the lines $l$, $g_1$, $g_2$, $g_3$, $g_4$ are concurrent.

\end{theorem}

\begin{theorem}[\textbf{Steiner lines}]

Consider any four lines $l_1$, $l_2$, $l_3$, $l_4$. Let $l$ be the Steiner line of complete quadrilateral formed by these four lines. Let $g_1$ be the Steiner line of the complete quadrilateral formed by four lines $l$, $l_2$, $l_3$, $l_4$. Similarly define the lines $g_2$, $g_3$, $g_4$. Then the Steiner line of the complete quadrilateral formed by $g_1$, $g_2$, $g_3$, $g_4$ is parallel to $l$.

\end{theorem}

\begin{theorem}[\textbf{Euler lines}]

Let given triangle $ABC$, let the Euler line of $ABC$ cut it's sides at $A_1$, $B_1$, $C_1$ respectively. So we get three new triangles $AB_1C_1$, $A_1BC_1$, $A_1B_1C$. Let three Euler lines of $AB_1C_1$, $A_1BC_1$, $A_1B_1C$ form a triangle $\triangle$. Then we have that :

\begin{enumerate}

\item Euler lines of $ABC$ and $\triangle$ are the same.

\item Let the Euler line of $\triangle$ and three side lines of $\triangle$ form triangles $\triangle$, $\triangle_1$, $\triangle_2$, $\triangle_3$. Then the Euler lines of $\triangle_1$, $\triangle_2$, $\triangle_3$ form a triangle $ABC$.

\end{enumerate}

\end{theorem}

\begin{theorem}[\textbf{Simson lines}]

Let given triangle $ABC$, let $P$ be any point on it's circumcircle. Let the Simson line of $P$ wrt $ABC$ intersects with $(ABC)$ at $Q$, $R$. Let the Simson lines of $Q$, $R$ wrt $ABC$ intersect at $N$. Then the middle of $PN$ lies on $QR$.

\end{theorem}

\section{On Three Pascal lines}

Consider some nice situation which includes three Pascal lines $l_1, l_2$, $l_3$ of some nice hexagons. Then we can predict that these lines are concurrent at the same point.

\begin{theorem}

Let given triangle $ABC$ with orthocenter $H$. Let $AA'$, $BB'$, $CC'$ be altitudes of $ABC$. Let $I_1$, $I_2$ be the incenters of the triangles $C'B'A$, $C'A'B$ respectively. Let $I_1'$, $I_2'$ be $C'$~-- excenter of $AC'B'$ and respectively $C'$~-- excenter of $C'BA'$. From theorem 6.6 we know that the points $I_1$, $I_2$, $I_1'$, $I_2'$, $A$, $B$ lie on the same conic, so let $l_C$ be the Pascal line of hexagon $AI_1I_1'BI_2I_2'$ in this order. Similarly define the lines $l_A$, $l_B$. Then $l_A$, $l_B$, $l_C$ are concurrent.

\end{theorem}

\begin{theorem}

Consider triangle $ABC$ and it's first Fermat point $F$. Let $O_A$ be the circumcenter of $BFC$. Similarly define the points $O_B$, $O_C$. Let $X_A$ be the reflection of $O_A$ wrt side $BC$. Similarly define the points $X_B$, $X_C$. It's well known that circles $(O_ABO_C)$, $(O_AO_BC)$ and $(AO_BO_C)$ intersect at the second Fermat point $F_2$. Let the lines $F_2X_A$, $F_2X_B$ intersect with $(X_AX_BC)$ second time at $T_{CA}$, $T_{CB}$ respectively. Similarly define the points $T_{AB}$, $T_{AC}$, $T_{BA}$, $T_{BC}$. Let $l_C$ be the Pascal line of hexagon $T_{CA}CT_{CB}O_AO_B$ in this order. Similarly define $l_A$, $l_B$. Then $l_A$, $l_B$, $l_C$ are concurrent.

\end{theorem}

\begin{theorem} Consider any points $A'$, $B'$, $C'$ on the sides $BC$, $CA$, $AB$ of $ABC$ respectively. Let the circles $(A'B'C)$, $(AB'C')$ and $(A'BC')$ intersect at $X$. Let a line through point $X$ and perpendicular to $AB$ intersects with $(AB'C')$ at $P_{AB}$ and a line through point $X$ and perpendicular to $AC$ intersects with $(AB'C')$ at $P_{AC}$. Similarly define the points $P_{BA}$, $P_{BC}$ on $(A'BC')$, and $P_{CA}$, $P_{CB}$ on $(A'B'C)$. Let $l_A$ be the Pascal line of hexagon $P_{AB}XP_{AC}B'AC'$ in this order. Similarly define $l_B$, $l_C$. Then we have that :

\begin{enumerate}

\item Lines $AB$, $l_A$, $XP_{AC}$ are concurrent.

\item Let $A'=l_B\cap l_C$, $B' = l_A\cap l_C$, $C' = l_A\cap l_B$. Then $XP_{BC}$, $BB'$, $CC'$ are concurrent.

\end{enumerate}

\end{theorem}

\subsection{New Pascal line} Consider any conic $\mathcal{C}$ and let conics $\mathcal{C}_1$, $\mathcal{C}_2$, $\mathcal{C}_3$ are tangent to it at six points. Let the common internal tangents to the conics $\mathcal{C}_1$, $\mathcal{C}_2$ meet at $X_{12}$. Let the common internal tangents to conics $\mathcal{C}_1$, $\mathcal{C}_3$ meet at $X_{13}$. Let the common internal tangents to the conics $\mathcal{C}_2$, $\mathcal{C}_3$ meet at $X_{23}$. Then $X_{12}$, $X_{23}$, $X_{13}$ are collinear. See picture below for more details.

\medskip

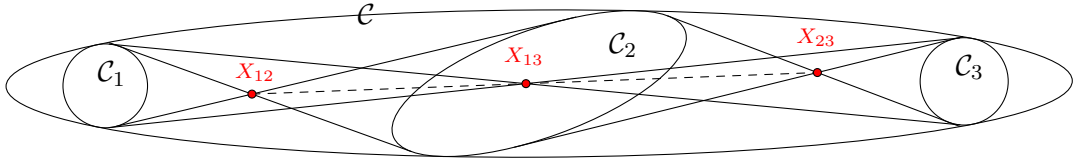
\definecolor{ffqqqq}{rgb}{1.0,0.0,0.0}
\begin{tikzpicture}[line cap=round,line join=round,>=triangle 45,x=1.0cm,y=1.0cm]
%\clip(-12.059535219428243,-10.703046786329853) rectangle (3.9182881984045768,7.913918927604378);
\draw [rotate around={0.3005284740075708:(-4.4782212554020795,0.028220915665245875)}] (-4.4782212554020795,0.028220915665245875) ellipse (7.0117514894687645cm and 0.9761705316932237cm);
\draw(-10.185053638503375,-0.0017129009385687926) circle (0.556042684842982cm);
\draw(1.1111860447961652,0.05753880724089326) circle (0.5791388509422064cm);
\draw [rotate around={19.990518003323153:(-4.4782212554020715,0.028220915665248227)}] (-4.4782212554020715,0.028220915665248227) ellipse (2.042931611061124cm and 0.7123716021331129cm);
\draw (-6.157168010203404,-0.8737789436285963)-- (-10.13250141524023,0.570927320190138);
\draw (-10.132079293368546,0.5518005913314532)-- (1.0560111488900652,-0.5189657803665344);
\draw (1.0499665189245055,0.6334328759051819)-- (-10.12627567576413,-0.5546402184619754);
\draw (1.04761934831359,-0.5381655820455087)-- (-2.788389514350013,0.9261663033436611);
\draw (-4.802315373508265,-0.815100518573936)-- (1.0929852475032855,0.6497565204647491);
\draw [dash pattern=on 3pt off 3pt] (-8.255650450335173,-0.11115341845106604)-- (-0.8194503426661027,0.1745568799268181);
\draw (-6.997637111239369,1.219462078415699) node[anchor=north west] {$\mathcal{C}$};
\draw (-10.417316367534905,0.4543812956512786) node[anchor=north west] {$\mathcal{C}_1$};
\draw (-3.6938791856657147,0.8369216870334889) node[anchor=north west] {$\mathcal{C}_2$};
\draw (0.865693156061668,0.539390271513992) node[anchor=north west] {$\mathcal{C}_3$};
\draw (-10.055477174092994,-0.5424470382355465)-- (-3.995719911646171,0.909656439414368);
\begin{scriptsize}
\draw [fill=ffqqqq] (-8.255650450335173,-0.11115341845106604) circle (1.5pt);
\draw[color=ffqqqq] (-8.214811083819137,0.17810212409746012) node {$X_{12}$};
\draw [fill=ffqqqq] (-0.8194503426661027,0.1745568799268181) circle (1.5pt);
\draw[color=ffqqqq] (-0.8344863611926102,0.6456514913423838) node {$X_{23}$};
\draw [fill=ffqqqq] (-4.651849261824445,0.027310191720221075) circle (1.5pt);
\draw[color=ffqqqq] (-4.6792104968017165,0.3906245637542436) node {$X_{13}$};
\end{scriptsize}
\end{tikzpicture}
\medskip
\begin{definition}

Name this line as $\text{Pasc}_{\mathcal{C}}(\mathcal{C}_1, \mathcal{C}_2, \mathcal{C}_3)$, where order of conics is important.

\end{definition}

\begin{theorem}

Consider three conics $\mathcal{C}_A$, $\mathcal{C}_B$, $\mathcal{C}_C$. Let given that there exists a conic $\mathcal{C}$, which is tangent to them at six points. Let two outer tangents to the conics $\mathcal{C}_A$, $\mathcal{C}_B $ intersect with two outer tangents to the conics $\mathcal{C}_A $, $\mathcal{C}_C $ at the points $P_A^{(1)}$, $P_A^{(2)}$, $P_A^{(3)}$, $P_A^{(4)}$. Similarly define $P_B^{(1)}$, $P_B^{(2)}$, $P_B^{(3)}$, $P_B^{(4)}$, $P_C^{(1)}$, $P_C^{(2)}$, $P_C^{(3)}$, $P_C^{(4)}$ (see picture below for more details). Let by definition union of the external tangents to conics $\mathcal{C}_A$, $\mathcal{C}_B$ form a conic $\mathcal{C}_{AB}$. Similarly define conics $\mathcal{C}_{BC}$, $\mathcal{C}_{CA}$. Note that we can look on the segments $XY$ as on the degenerated conics. Let the lines $\text{Pasc}_{\mathcal{C}_{AB}}(\mathcal{C}_A, P_{A}^{(1)}P_{A}^{(4)}, \mathcal{C}_B)$ and $\text{Pasc}_{\mathcal{C}_{BC}}(\mathcal{C}_C, P_{C}^{(1)}P_{C}^{(2)}, \mathcal{C}_B)$ meet at $W_B$. Similarly define $W_A$, $W_C$. Then $P_{A}^{(1)}W_A$, $P_{B}^{(1)}W_B$, $P_{C}^{(1)}W_C$ are concurrent.

\end{theorem}

\newpage

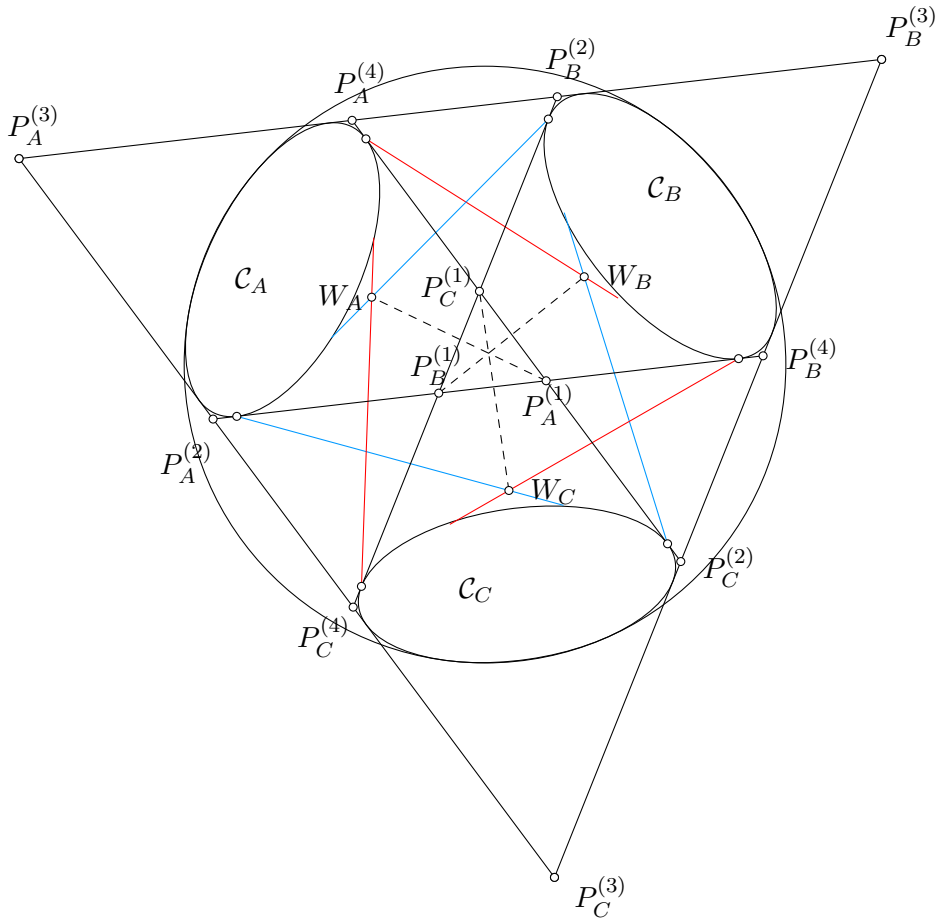
\definecolor{ffqqqq}{rgb}{1.0,0.0,0.0}
\definecolor{qqccqq}{rgb}{0.0,0.6,1.0}
\definecolor{ffffff}{rgb}{1.0,1.0,1.0}
\begin{tikzpicture}[line cap=round,line join=round,>=triangle 45,x=1.0cm,y=1.0cm]
%\clip(-13.28355449956431,-9.469563850599604) rectangle (0.8041230571397413,6.945050531069062);
\draw(-5.927380554700882,-0.7956155059193679) circle (3.956404344821457cm);
\draw [rotate around={-51.47741441083869:(-3.6237416740074337,1.038264699620901)}] (-3.6237416740074337,1.038264699620901) ellipse (2.102833721950116cm and 1.0042143224650857cm);
\draw [rotate around={64.70652146073907:(-8.589562353057218,0.4624213669167215)}] (-8.589562353057218,0.4624213669167215) ellipse (2.1028337219500925cm and 1.0042143224650744cm);
\draw [rotate around={8.207378488465768:(-5.507039478345662,-3.709921955221418)}] (-5.507039478345662,-3.709921955221418) ellipse (2.1028337219501028cm and 1.0042143224650792cm);
\draw [dash pattern=on 3pt off 3pt] (-5.123694065825692,-1.0104834813962966)-- (-7.417872926140182,0.09780563653314696);
\draw [dash pattern=on 3pt off 3pt] (-6.535834784234573,-1.1742372433810018)-- (-4.622566169429861,0.36892907894155635);
\draw [dash pattern=on 3pt off 3pt] (-5.613063098868942,-2.465201944792031)-- (-6.000277484387625,0.17601443650678822);
\draw [color=qqccqq] (-9.191971609609855,-1.4822464890823759)-- (-4.899286421260857,-2.6612423513314907);
\draw [color=ffqqqq] (-7.552245361125632,-3.7368203494185392)-- (-7.390424527599698,0.8811087174626273);
\draw [color=ffqqqq] (-7.493646738091285,2.1973617180772216)-- (-4.181240525504717,0.0878731515502491);
\draw [color=qqccqq] (-5.094482956019033,2.459711407703674)-- (-7.9447301764669085,-0.4377855453745729);
\draw [color=qqccqq] (-3.5266826674854332,-3.1721154024139064)-- (-4.88436522045534,1.2148603145272876);
\draw [color=ffqqqq] (-2.5922903333563,-0.7169384551179991)-- (-6.385380304968094,-2.912178291137984);
\draw (-12.055430640288748,1.9353262152297415)-- (-0.711648563780299,3.250766642622808);
\draw (-0.711648563780299,3.250766642622808)-- (-5.01380147230369,-7.5958583469496315);
\draw (-5.01380147230369,-7.5958583469496315)-- (-12.055430640288748,1.9353262152297415);
\draw (-7.661447520642244,-4.012141775112518)-- (-4.976900122983826,2.756162259698645);
\draw (-7.675309504667428,2.4432510307643285)-- (-3.3526314360490708,-3.407702135330326);
\draw (-9.503815201447011,-1.5184082969308847)-- (-2.2705832250310456,-0.6796328604568388);
\draw (-9.348544202588934,0.6115577673658549) node[anchor=north west] {$\mathcal{C}_A$};
\draw (-3.914482363908655,1.829537145001087) node[anchor=north west] {$\mathcal{C}_B$};
\draw (-6.401545149009974,-3.492095828051312) node[anchor=north west] {$\mathcal{C}_C$};
\draw (-12.346647285998746,2.7851825028379613) node[anchor=north west] {$P_A^{(3)}$};
\draw (-8.053908780207426,3.347326830977299) node[anchor=north west] {$P_A^{(4)}$};
\draw (-10.336555445985349,-1.6932339780054304) node[anchor=north west] {$P_A^{(2)}$};
\draw (-5.277256492731295,3.6658752835895907) node[anchor=north west] {$P_B^{(2)}$};
\draw (-0.7971365442268638,4.09685260182975) node[anchor=north west] {$P_B^{(3)}$};
\draw (-2.1088066432186556,-0.30661130192839703) node[anchor=north west] {$P_B^{(4)}$};
\draw (-6.912585447318464,0.6490340559084774) node[anchor=north west] {$P_C^{(1)}$};
\draw (-7.048862860200728,-0.40030202328495335) node[anchor=north west] {$P_B^{(1)}$};
\draw (-5.58388067171639,-0.9249700628816687) node[anchor=north west] {$P_A^{(1)}$};
\draw (-3.1990259462767683,-3.079856654082464) node[anchor=north west] {$P_C^{(2)}$};
\draw (-4.885458930694786,-7.4458442692979885) node[anchor=north west] {$P_C^{(3)}$};
\draw (-8.53087972529535,-3.998025723376716) node[anchor=north west] {$P_C^{(4)}$};
\draw (-5.464159775438935,-2.1846707115681614) node[anchor=north west] {$W_C$};
\draw (-8.260440798003908,0.39112849060047655) node[anchor=north west] {$W_A$};
\draw (-4.464342979860318,0.662265248723491) node[anchor=north west] {$W_B$};
\begin{scriptsize}
\draw [fill=ffffff] (-4.622566169429861,0.36892907894155635) circle (1.5pt);
\draw [fill=ffffff] (-7.417872926140182,0.09780563653314696) circle (1.5pt);
\draw [fill=ffffff] (-5.613063098868942,-2.465201944792031) circle (1.5pt);
\draw [fill=ffffff] (-9.191971609609855,-1.4822464890823759) circle (1.5pt);
\draw [fill=ffffff] (-7.552245361125632,-3.7368203494185392) circle (1.5pt);
\draw [fill=ffffff] (-7.493646738091285,2.1973617180772216) circle (1.5pt);
\draw [fill=ffffff] (-5.094482956019033,2.459711407703674) circle (1.5pt);
\draw [fill=ffffff] (-3.5266826674854332,-3.1721154024139064) circle (1.5pt);
\draw [fill=ffffff] (-2.5922903333563,-0.7169384551179991) circle (1.5pt);
\draw [fill=ffffff] (-12.055430640288748,1.9353262152297415) circle (1.5pt);
\draw [fill=ffffff] (-0.711648563780299,3.250766642622808) circle (1.5pt);
\draw [fill=ffffff] (-0.711648563780299,3.250766642622808) circle (1.5pt);
\draw [fill=ffffff] (-5.01380147230369,-7.5958583469496315) circle (1.5pt);
\draw [fill=ffffff] (-5.01380147230369,-7.5958583469496315) circle (1.5pt);
\draw [fill=ffffff] (-12.055430640288748,1.9353262152297415) circle (1.5pt);
\draw [fill=ffffff] (-7.661447520642244,-4.012141775112518) circle (1.5pt);
\draw [fill=ffffff] (-4.976900122983826,2.756162259698645) circle (1.5pt);
\draw [fill=ffffff] (-7.675309504667428,2.4432510307643285) circle (1.5pt);
\draw [fill=ffffff] (-3.3526314360490708,-3.407702135330326) circle (1.5pt);
\draw [fill=ffffff] (-9.503815201447011,-1.5184082969308847) circle (1.5pt);
\draw [fill=ffffff] (-2.2705832250310456,-0.6796328604568388) circle (1.5pt);
\draw [fill=ffffff] (-6.535834784234573,-1.1742372433810018) circle (1.5pt);
\draw [fill=ffffff] (-6.000277484387625,0.17601443650678822) circle (1.5pt);
\draw [fill=ffffff] (-5.123694065825692,-1.0104834813962966) circle (1.5pt);
\end{scriptsize}
\end{tikzpicture}

\section{Facts related to the set of confocal conics}

\begin{theorem}

Consider any set of points $A_1$, $A_2$, $\ldots$, $A_n$. Let given a set of conics $\mathcal{K}_1$, $\mathcal{K}_2$, $\ldots$, $\mathcal{K}_n$, where the conic $\mathcal{K}_i$ has foci $A_{i}$ and $A_{i+1}$, for every $1\leq i\leq n$ ($A_{n+1} = A_1$). For each integer $i$, consider the line $l_i$ passing through two intersection points of the conics $\mathcal{K}_i$, $\mathcal{K}_{i + 1}$. Then in the case when the lines $l_1$, $l_2$,$\ldots$ , $l_{n-1}$ pass through the same point $P$ we get that the line $l_n$ also passes through $P$.

\end{theorem}

\begin{commentary}

Previous theorem can be seen as $n$~- conic analog of \cite[Problem 11.18]{Ak}.

\end{commentary}

By observing theorem 13.1 one can note that $P$ can be seen as the incenter for the set of conics $\mathcal{K}_1$, $\mathcal{K}_2$, $\ldots$, $\mathcal{K}_n$. So if we apply this analogy to \cite[Problem 5.5.10]{Ak},  then we'll obtain the following theorem.

\begin{theorem}

Consider cyclic quadrilateral $ABCD$. Let $\mathcal{K}_{AB}$ be an ellipse with foci $A$ and $B$. Let $\mathcal{K}_{CD}$ be an ellipse with foci $C$ and $D$. Let the segments $AC$, $AD$, $BC$, $BD$ intersect with the ellipse $\mathcal{K}_{AB}$ at $T_{AC}$, $T_{AD}$, $T_{BC}$, $T_{BD}$ respectively. Let the segments $CA$, $CB$, $DA$, $DB$ intersect with the ellipse $\mathcal{K}_{CD}$ at $T_{CA}$, $T_{CD}$, $T_{DA}$, $T_{DB}$ respectively. Consider the intersection point $I_C$ of the tangent lines through $T_{AC}$ and $T_{BC}$ to $\mathcal{K}_{AB}$. Consider the intersection point $I_D$ of the tangent lines through $T_{AD}$ and $T_{BD}$ to $\mathcal{K}_{AB}$. Similarly define the points $I_A$, $I_B$ (similarly for ellipse $\mathcal{K}_{CD}$). Then the points $I_A$, $I_B$, $I_C$, $I_D$ lie on the same circle.

\end{theorem}

We can look at the set of three confocal conics (same situation as \cite[Problem 11.18]{Ak}) as on the analog of the situation where given three circles $\omega_A$, $\omega_B$, $\omega_C$ with centers at $A$, $B$ and $C$. So \cite[Problem 11.18]{Ak} can be seen as an analog of radical center theorem. Next two pictures describes conic analogs of \cite[Theorem 6.3.7]{Ak} and \cite[Theorem 10.11]{Ak}.

\bigskip

\definecolor{ffffff}{rgb}{1.0,1.0,1.0}
\begin{tikzpicture}[line cap=round,line join=round,>=triangle 45,x=1.0cm,y=1.0cm]
%\clip(-7.274272441060294,-9.333962096203539) rectangle (23.833748032781706,9.7614110251098);
\draw [rotate around={-66.1707356271919:(-2.1245186707387655,2.1147289690092292)}] (-2.1245186707387655,2.1147289690092292) ellipse (1.9599950257357077cm and 0.8424823523000249cm);
\draw [rotate around={63.8345250658968:(-3.5233589298036194,2.341649347394238)}] (-3.5233589298036194,2.341649347394238) ellipse (2.3343746602986277cm and 1.7447701756644716cm);
\draw [rotate around={-9.214275880875809:(-2.808381780666485,0.7228190094083604)}] (-2.808381780666485,0.7228190094083604) ellipse (2.1246968422389774cm and 1.5830634045777392cm);
\draw [dash pattern=on 3pt off 3pt] (-6.309756829868154,-0.2867912419496193)-- (-3.5477470997927134,1.3375859844858176);
\draw [dash pattern=on 3pt off 3pt] (-3.5477470997927134,1.3375859844858176)-- (0.25212733808678844,-0.15820389265697415);
\draw [dash pattern=on 3pt off 3pt] (-3.5477470997927134,1.3375859844858176)-- (-1.9992752530123383,6.575977132862175);
\draw (-6.309756829868154,-0.2867912419496193)-- (-1.4095415216016314,0.495898631023351);
\draw (-1.4095415216016314,0.495898631023351)-- (-1.9992752530123383,6.575977132862175);
\draw (-2.839495819875898,3.7335593069951054)-- (0.25212733808678844,-0.15820389265697415);
\draw (-4.207222039731341,0.9497393877933701)-- (0.25212733808678844,-0.15820389265697415);
\draw (-6.309756829868154,-0.2867912419496193)-- (-2.839495819875898,3.7335593069951054);
\draw (-4.207222039731341,0.9497393877933701)-- (-1.9992752530123383,6.575977132862175);
\draw [rotate around={68.19859051364821:(3.6404447770974184,2.0328440008127457)}] (3.6404447770974184,2.0328440008127457) ellipse (1.2617404435303345cm and 1.0526824170408495cm);
\draw [rotate around={-66.06525220653162:(4.225093515501749,1.9435238524373069)}] (4.225093515501749,1.9435238524373069) ellipse (1.3904850022530333cm and 1.1342427325008433cm);
\draw [rotate around={-8.686240940738992:(3.9667570132968324,1.2976825969250183)}] (3.9667570132968324,1.2976825969250183) ellipse (1.1610084026014185cm and 0.9990737083356319cm);
\draw (1.3817593984566303,0.34031303529374074)-- (5.898312579549596,0.2617144881301778);
\draw [dash pattern=on 3pt off 3pt] (5.898312579549596,0.2617144881301778)-- (3.956357578684052,6.968743752474139);
\draw (3.8595816544703654,4.4439674717053705)-- (7.484444689637288,-0.6823441563935639);
\draw (2.173398509998168,0.742081063360782)-- (7.484444689637288,-0.6823441563935639);
\draw (2.173398509998168,0.742081063360782)-- (3.956357578684052,6.968743752474139);
\draw (1.3817593984566303,0.34031303529374074)-- (3.8595816544703654,4.4439674717053705);
\draw (3.720758942899589,0.8222372505401745)-- (3.956357578684052,6.968743752474139);
\draw (3.057074428122954,1.952806561070234)-- (7.484444689637288,-0.6823441563935639);
\draw (4.7158641976714035,2.0324183027317453)-- (1.3817593984566303,0.34031303529374074);
\begin{scriptsize}
\draw [fill=ffffff] (-1.4095415216016314,0.495898631023351) circle (1.5pt);
\draw [fill=ffffff] (-2.839495819875898,3.7335593069951054) circle (1.5pt);
\draw [fill=ffffff] (-4.207222039731341,0.9497393877933701) circle (1.5pt);
\draw [fill=ffffff] (-1.9992752530123383,6.575977132862175) circle (1.5pt);
\draw [fill=ffffff] (-6.309756829868154,-0.2867912419496193) circle (1.5pt);
\draw [fill=ffffff] (0.25212733808678844,-0.15820389265697415) circle (1.5pt);
\draw [fill=ffffff] (-3.5477470997927134,1.3375859844858176) circle (1.5pt);
\draw [fill=ffffff] (3.3821082748925018,1.3870027453004565) circle (1.5pt);
\draw [fill=ffffff] (3.898781279302334,2.678685256325034) circle (1.5pt);
\draw [fill=ffffff] (4.551405751701162,1.2083624485495799) circle (1.5pt);
\draw [fill=ffffff] (4.7158641976714035,2.0324183027317453) circle (1.5pt);
\draw [fill=ffffff] (2.8478783090062776,1.0843890858773324) circle (1.5pt);
\draw [fill=ffffff] (3.057074428122954,1.952806561070234) circle (1.5pt);
\draw [fill=ffffff] (4.9958929673849495,0.7988305520655441) circle (1.5pt);
\draw [fill=ffffff] (3.720758942899589,0.8222372505401745) circle (1.5pt);
\draw [fill=ffffff] (3.8144695012973817,3.2670415448181815) circle (1.5pt);
\draw [fill=ffffff] (3.956357578684052,6.968743752474139) circle (1.5pt);
\draw [fill=ffffff] (2.173398509998168,0.742081063360782) circle (1.5pt);
\draw [fill=ffffff] (7.484444689637288,-0.6823441563935639) circle (1.5pt);
\draw [fill=ffffff] (3.8595816544703654,4.4439674717053705) circle (1.5pt);
\draw [fill=ffffff] (1.3817593984566303,0.34031303529374074) circle (1.5pt);
\draw [fill=ffffff] (5.898312579549596,0.2617144881301778) circle (1.5pt);
\end{scriptsize}
\end{tikzpicture}

\bigskip

\begin{definition}

Consider any two conics $\mathcal{K}_1$ and $\mathcal{K}_2$ which share same focus $F$. Let the conics $\mathcal{K}_1$ and $\mathcal{K}_2$ intersect at the points $A$, $B$. Consider the intersection point $X$ of the tangents to the conic $\mathcal{K}_1$ from $A$, $B$. Similarly let $Y$ be the intersection point of the tangents to the conic $\mathcal{K}_2$ from  $A$, $B$. Then let by definition $\mathcal{L}_F(\mathcal{K}_1 , \mathcal{K}_2) = XY$.
 
\end{definition}

\begin{theorem}

Consider points $A$, $B$, $C$. Let given set of conics $\mathcal{K}_{AB}$, $\mathcal{K}_{BC}$, $\mathcal{K}_{CA}$, where conic $\mathcal{K}_{AB}$ has foci $A$ and $B$, similarly for $\mathcal{K}_{BC}$, $\mathcal{K}_{CA}$. Let $X_C$ be the intersection point of the lines $\mathcal{L}_{A}(\mathcal{K}_{AB} , \mathcal{K}_{CA})$ and $\mathcal{L}_{B}(\mathcal{K}_{AB} , \mathcal{K}_{BC})$. Similarly define the points $X_A$ and $X_B$. Then $AX_A$, $BX_B$ and $CX_C$ are concurrent.

\end{theorem}

\begin{theorem}

Consider cyclic quadrilateral $ABCD$. Let $P= AC\cap BD$. Let conic $\mathcal{K}_{AB}$ has foci at points $A$, $B$ and goes through the incenter of $ABP$. Like the same let conic $\mathcal{K}_{BC}$ has foci at points $B$, $C$ and goes through the incenter of $BCP$. Similarly define conics $\mathcal{K}_{CD}$ and $\mathcal{K}_{DA}$. Consider $P_{AB}$ as the intersection point of the lines $\mathcal{L}_A(\mathcal{K}_{AB} , \mathcal{K}_{DA})$ and $\mathcal{L}_B(\mathcal{K}_{AB} , \mathcal{K}_{BC})$. Similarly define the points $P_{BC}$, $P_{CD}$, $P_{DA}$. Then we have that :

\begin{enumerate}

\item Lines $P_{AB}P_{CD}$, $P_{BC}P_{DA}$ intersect at $P$.

\item Define point $A^*$ as the second intersection point of the line $\mathcal{L}_A(\mathcal{K}_{AB} , \mathcal{K}_{DA})$ with $(ABCD)$. Similarly define the points $B^*$, $C^*$, $D^*$. Then the lines $A^*C^*$ and $B^*D^*$ intersect at $P$.

\end{enumerate}

\end{theorem}

\newpage

\definecolor{ffffff}{rgb}{1.0,1.0,1.0}
\begin{tikzpicture}[line cap=round,line join=round,>=triangle 45,x=1.0cm,y=1.0cm]
%\clip(-13.564408706003515,-1.1937811941018512) rectangle (14.237932299703608,15.736553395784647);
\draw(-1.9449837908720675,9.527425330318346) circle (3.9774050150527cm);
\draw (-0.7384668608292344,5.737429611153694)-- (-2.0215535384729044,13.50409324800002);
\draw (-5.653797867363089,10.964244316756238)-- (1.8559109979214456,10.699154590125698);
\draw (1.8559109979214456,10.699154590125698)-- (-2.0215535384729044,13.50409324800002);
\draw (-2.0215535384729044,13.50409324800002)-- (-5.653797867363089,10.964244316756238);
\draw (-5.653797867363089,10.964244316756238)-- (-0.7384668608292344,5.737429611153694);
\draw (-0.7384668608292344,5.737429611153694)-- (1.8559109979214456,10.699154590125698);
\draw [rotate around={62.395912549730156:(0.5587220685461054,8.218292100639694)}] (0.5587220685461054,8.218292100639694) ellipse (3.0756734017851968cm and 1.273730373437171cm);
\draw [rotate around={-35.881791064018394:(-0.08282127027571821,12.101623919062845)}] (-0.08282127027571821,12.101623919062845) ellipse (2.537916787245108cm and 0.8458242090544579cm);
\draw [rotate around={34.9631919592365:(-3.837675702918005,12.234168782378124)}] (-3.837675702918005,12.234168782378124) ellipse (2.4346702009011882cm and 1.0082713536457595cm);
\draw [rotate around={-46.759104427014925:(-3.1961323640961674,8.350836963954976)}] (-3.1961323640961674,8.350836963954976) ellipse (3.805014656231591cm and 1.2681162472396246cm);
\draw [dash pattern=on 3pt off 3pt] (-3.459445247075839,5.849632901232724)-- (-0.6449055271618617,13.2863543216661);
\draw (-5.377995114614327,14.034968227723823)-- (0.542311933802695,13.098576870094348);
\draw (0.542311933802695,13.098576870094348)-- (4.694942127371137,5.513375510649789);
\draw (4.694942127371137,5.513375510649789)-- (-6.10338374815116,5.95865934557596);
\draw (-6.10338374815116,5.95865934557596)-- (-5.377995114614327,14.034968227723823);
\draw [dash pattern=on 3pt off 3pt] (-5.850473881610498,8.774499168659519)-- (1.0900296731508101,12.098114717810779);
\draw [dash pattern=on 3pt off 3pt] (4.694942127371137,5.513375510649789)-- (-5.377995114614327,14.034968227723823);
\draw [dash pattern=on 3pt off 3pt] (-6.10338374815116,5.95865934557596)-- (0.542311933802695,13.098576870094348);
\draw(-2.72633644771811,11.829772883141644) circle (0.968263716042755cm);
\draw(-0.8637345068813359,11.633041766565189) circle (0.8373632344715588cm);
\draw(-0.12776239631046318,9.545226149946167) circle (1.2231894719227405cm);
\draw(-2.649354297811856,9.601975740657382) circle (1.2554309880318006cm);
\begin{scriptsize}
\draw [fill=ffffff] (-0.7384668608292344,5.737429611153694) circle (1.5pt);
\draw [fill=ffffff] (-5.653797867363089,10.964244316756238) circle (1.5pt);
\draw [fill=ffffff] (-2.0215535384729044,13.50409324800002) circle (1.5pt);
\draw [fill=ffffff] (1.8559109979214456,10.699154590125698) circle (1.5pt);
\draw [fill=ffffff] (-1.5781919282523986,10.820377053239929) circle (1.5pt);
\draw [fill=ffffff] (-0.12776239631046318,9.545226149946167) circle (1.5pt);
\draw [fill=ffffff] (-0.8637345068813359,11.633041766565189) circle (1.5pt);
\draw [fill=ffffff] (-2.72633644771811,11.829772883141644) circle (1.5pt);
\draw [fill=ffffff] (-2.649354297811856,9.601975740657382) circle (1.5pt);
\draw [fill=ffffff] (-3.459445247075839,5.849632901232724) circle (1.5pt);
\draw [fill=ffffff] (-0.6449055271618617,13.2863543216661) circle (1.5pt);
\draw [fill=ffffff] (-5.377995114614327,14.034968227723823) circle (1.5pt);
\draw [fill=ffffff] (0.542311933802695,13.098576870094348) circle (1.5pt);
\draw [fill=ffffff] (4.694942127371137,5.513375510649789) circle (1.5pt);
\draw [fill=ffffff] (-6.10338374815116,5.95865934557596) circle (1.5pt);
\draw [fill=ffffff] (-5.850473881610498,8.774499168659519) circle (1.5pt);
\draw [fill=ffffff] (1.0900296731508101,12.098114717810779) circle (1.5pt);
\end{scriptsize}
\end{tikzpicture}

\end{document}